# ASYMPTOTIC THEOREMS OF SEQUENTIAL ESTIMATION-ADJUSTED URN MODELS

By Li-X. Zhang[1], Feifang Hu[2] and Siu Hung Cheung[3]

*Zhejiang University, University of Virginia and Chinese University of Hong Kong*

The Generalized Pólya Urn (GPU) is a popular urn model which is widely used in many disciplines. In particular, it is extensively used in treatment allocation schemes in clinical trials. In this paper, we propose a sequential estimation-adjusted urn model (a nonhomogeneous GPU) which has a wide spectrum of applications. Because the proposed urn model depends on sequential estimations of unknown parameters, the derivation of asymptotic properties is mathematically intricate and the corresponding results are unavailable in the literature. We overcome these hurdles and establish the strong consistency and asymptotic normality for both the patient allocation and the estimators of unknown parameters, under some widely satisfied conditions. These properties are important for statistical inferences and they are also useful for the understanding of the urn limiting process. A superior feature of our proposed model is its capability to yield limiting treatment proportions according to any desired allocation target. The applicability of our model is illustrated with a number of examples.

## 1. Introduction.

1.1. *Brief history.* Urn models have long been recognized as valuable mathematical apparatus in many areas including physical sciences, biological sciences and engineering [19, 20]. Among various urn models, the Pólya urn

Received March 2005.
[1]Supported by National Natural Science Foundation of China Grant 10471126.
[2]Supported in part by NSF Grants DMS-02-04232, DMS-03-49048 and the Research Grants Council of the Hong Kong Special Administrative Region Project CUHK400204.
[3]Supported by the Research Grants Council of the Hong Kong Special Administrative Region Project CUHK400204.

(also known as the Pólya–Eggenberger urn) model is particularly popular. It was originally proposed to model the problem of contagious diseases [10]. Since then, there have been numerous generalizations and extensions. For example, in information science, Martin and Ho [21] apply the Pólya urn process to study human decision-making. In population genetics, Benaïm, Schreiber and Tarrés [9] make use of a class of generalized Pólya urn models to scrutinize the evolutionary processes. In economics, the model is employed to capture the mechanism of reinforcement learning by Erev and Roth [13] and Beggs [8]. In business, a coupled Pólya urn model is utilized by Windrum [32] to examine the recent browser war between Netscape and Microsoft.

Urn models are also extensively applied in clinical studies. Applications are mostly found in the area of adaptive design in which urn models are utilized to provide response-adaptive allocation schemes. In most clinical trials, patients accrue sequentially. Assume the availability of several treatments. Adaptive designs are valuable and ethical randomization schemes that formulate treatment allocation as a function of previous responses, with one of the important objectives to assign more patients to the better treatment.

Early works in adaptive design can be traced to Thompson [29] and Robbins [24]. Athreya and Karlin [3] successfully embed the urn schemes into a continuous-time branching process and generate important limit theorems. Among many useful allocation schemes generated under the paradigm of urn models, the most influential family of models is the Generalized Pólya Urn (GPU) [also named as the Generalized Friedman's Urn (GFU) in literature] [2, 3, 30].

A general description of the GPU model is as follows. To compare $K$ treatments in a clinical study, consider an urn containing particles of $K$ types, representing the $K$ treatments. Patients arrive sequentially to be allocated to the treatments. At the beginning, the urn contains $\mathbf{Y}_0 = (Y_{0,1}, \ldots, Y_{0,K})$ particles, where $Y_{0,k} > 0$ denotes the number of particles of type $k$, $k = 1, \ldots, K$. At stage $i$, $i = 1, 2, \ldots$, a particle is drawn from the urn and replaced. If the particle is of type $k$, then treatment $k$ is assigned to the $i$th patient, $k = 1, \ldots, K$, $i = 1, 2, \ldots$. After observing a random variable $\xi_{i,k}$ which represents the response of the $i$th patient on treatment $k$, additional $d_i(k, g, \xi_{i,k})$ particles of type $g$, $g = 1, \ldots, K$, are added to the urn, where $d_i(k, g, \xi_{i,k})$ is a function of $\xi_{i,k}$. After $n$ splits and generations, the urn composition is denoted by the row vector $\mathbf{Y}_n = (Y_{n,1}, \ldots, Y_{n,K})$, where $Y_{n,k}$ stands for the number of particles of type $k$ in the urn after the $n$th split. The urn composition in stage $n$ can be characterized by the following recursive formula:

$$(1.1) \qquad \mathbf{Y}_n = \mathbf{Y}_{n-1} + \mathbf{X}_n \mathbf{D}_n,$$

where the matrix $\mathbf{D}_n = (d_n(k, g, \xi_{n,k}), k, g = 1, \ldots, K)$ and $\mathbf{X}_n = (X_{n,1}, \ldots, X_{n,K})$ is the observed result of the $n$th draw, distributed according to the urn composition at the previous stage. That is, if the $n$th draw is a type-$k$ particle, all elements in $\mathbf{X}_n$ are 0 except the $k$th component $X_{n,k}$ which is 1. Moreover, define $\mathbf{H}_i = (\mathsf{E}[d_i(k, g, \xi_{i,k})|\mathcal{F}_{i-1}], k, q = 1, \ldots, K)$, where the sigma field $\mathcal{F}_i$ is generated by $\{\mathbf{Y}_0, \mathbf{Y}_1, \ldots, \mathbf{Y}_i, \mathbf{X}_1, \ldots, \mathbf{X}_i\}$. The matrices $\mathbf{D}_i$ and $\mathbf{H}_i$ are called the *addition rules* and the *generating matrices*, respectively. These matrices play important roles in the mechanism of the GPU model.

A GFU model is said to be *homogeneous* if $\mathbf{H}_i = \mathbf{H}$ for all $i = 1, \ldots, n$. For instance, if the addition rule $\mathbf{D}_i$ is merely a function of the $i$th patient's observed outcome, the addition rules are independent and identically distributed (i.i.d.). Consequently $\mathbf{H}_i = \mathbf{H} = \mathsf{E}\mathbf{D}_i$ are identical and nonrandom. However, in clinical studies, the entire history of all previous trials provides more information of the efficacies of the treatments. For applications of such models, one can refer to Andersen, Faries and Tamura [1] and Bai, Hu and Shen [6]. If the addition rules rely on all the previous outcomes, the GPU model is no longer homogeneous. Under these circumstances, we will only assume that $\mathbf{H}_n$ converges to a limit $\mathbf{H}$.

In clinical investigations, a crucial element of the treatment allocation process is the number of patients being assigned to a particular treatment. Here, we define $\mathbf{N}_n = (N_{n,1}, \ldots, N_{n,K})$. The quantity $N_{n,k}$ is the number of times a type-$k$ particle is selected in the first $n$ stages. In clinical trials, $N_{n,k}$ represents the number of patients assigned to treatment $k$ among the first $n$ patients. Obviously,

$$\mathbf{N}_n = \sum_{i=1}^{n} \mathbf{X}_i. \tag{1.2}$$

Owing to its clinical applications, the asymptotic behavior of $\mathbf{N}_n$ is of immense importance [15]. Athreya and Karlin [2, 3] obtained the asymptotic convergence of $\mathbf{Y}_n$ and $\mathbf{N}_n$ and conjecture that $\mathbf{N}_n$ is asymptotically normal. This conjecture has not been solved for almost three decades due to its mathematical complexity. A limited success was reported by Smythe [27] who derived the asymptotic normality of $\mathbf{Y}_n$ and $\mathbf{N}_n$ under the assumptions that the eigenvalues of the generating matrix $\mathbf{H}$ are simple and the GPU is homogeneous. Bai, Hu and Zhang [7] considered functional limit theorems of $\mathbf{Y}_n$ and $\mathbf{N}_n$ for $K = 2$. Further, Janson [17] established functional limit theorems of $\mathbf{Y}_n$ and $\mathbf{N}_n$ for the general homogeneous GPU model and the asymptotic normality of $\mathbf{N}_n$ follows.

1.2. *Objectives and organization of the paper.* Recently, Bai and Hu [5] established the consistency and the asymptotic normality of the GPU non-homogeneous model. However, their theorems rely on very stringent assumptions which are not valid in many practical applications (illustrated in Section 2). In this paper, we derive a sequential estimation-adjusted urn model (SEU) and successfully produce important asymptotic theorems related to $\mathbf{Y}_n$ and $\mathbf{N}_n$. Our allocation scheme is suitable for many popular clinical scenarios. Because the proposed urn model depends on sequential estimators of unknown parameters, we cannot use the techniques of Bai and Hu [5] to show the asymptotic properties. Instead, we show the law of the iterated logarithm of the estimators first and then prove the strong consistency and asymptotic normality of the SEUs by using the results of the law of the iterated logarithm. Also we obtain the asymptotic properties of the estimators of unknown parameters.

Another important advantage of our proposed SEU model is that it can be used to satisfy a pre-specified treatment allocation target. Different targets might be of interest in clinical studies [25]. In fact, there is a growing interest in target-based designs which are derived with a pre-specified allocation target (see, e.g., [11, 12, 16, 22, 23]). Based on some optimality consideration, a suitable allocation target, say $\mathbf{v}$, is used. It will be shown later that we are able to define an SEU model to ensure that the allocation proportion converges to $\mathbf{v}$ almost surely. Illustrative examples are given in Section 4.

To summarize, the major contributions of this paper are:

(a) to propose a general nonhomogeneous GPU model (SEU) with the addition rules utilizing sequential estimations of unknown parameters;
(b) to show strong consistency and asymptotic normality of $\mathbf{N}_n$ under widely satisfied conditions and obtain a general and explicit asymptotic variance of $\mathbf{N}_n$;
(c) to state the framework of the targeting procedures of SEU when the desired allocation proportion is pre-specified based on some optimality consideration; and
(d) to illustrate with examples the applications of the SEU model in clinical trials.

The SEU model will be described in Section 2. Asymptotic properties are presented in Section 3. Then in Section 4, selected examples are provided to illustrate the application of the SEU model in clinical trials. Concluding remarks are also given. Finally, technical proofs are stated in Section 5.

## 2. Sequential estimation-adjusted urn model.

2.1. *The SEU process.* Before the outline of the SEU model, we adopt the following assumption on patient responses. Let the responses $\{\{\xi_{j,k}, j = 1, 2, \ldots\}, k = 1, \ldots, K\}$ be $K$ independent sequences of i.i.d. random variables, and write $\boldsymbol{\xi}_j = (\xi_{j,1}, \ldots, \xi_{j,K})$. In clinical trial, $\xi_{i,k}$ can only be observed when the $i$th patient is assigned to treatment $k$, that is, when $X_{i,k} = 1$. Now, a useful nonhomogeneous GPU model with addition rules utilizing the entire history of the previous trials is given as follows.

*The SEU process.* Let $\boldsymbol{\Theta} = (\theta_1, \ldots, \theta_K)$ be the parameter from the distribution of the response $\boldsymbol{\xi}_n = (\xi_{n,1}, \ldots, \xi_{n,K})$. If the value of $\boldsymbol{\Theta}$ is known, it is natural to consider the adaptive design with addition rules $\mathbf{D}(\boldsymbol{\Theta}, \boldsymbol{\xi}_n)$ and generating matrices $\mathbf{H} = \mathbf{H}(\boldsymbol{\Theta}) = \mathsf{E}[\mathbf{D}(\boldsymbol{\Theta}, \boldsymbol{\xi}_n)]$, whose values depend on the parameter $\boldsymbol{\Theta}$. Here, the $d \times d$ matrix $\mathbf{H}(\mathbf{x})$ is a function continuous at point $\boldsymbol{\Theta}$, and also, as a function of $\mathbf{x}$, the matrix $\mathbf{D}(\mathbf{x}, \mathbf{y}) = \{D_{ij}(\mathbf{x}, \mathbf{y})\}$ is continuous at point $\boldsymbol{\Theta}$ for any fixed possible value $\mathbf{y}$ of $\boldsymbol{\xi}_n$. However, $\boldsymbol{\Theta}$ is usually unknown in practice; an estimate of $\boldsymbol{\Theta}$ is therefore required. Hence, without loss of generality, we assume $\boldsymbol{\Theta} = \mathsf{E}[\boldsymbol{\xi}_n]$, with its estimate $\widehat{\boldsymbol{\Theta}}_n = (\hat{\theta}_{n,1}, \ldots, \hat{\theta}_{n,K})$, where $\hat{\theta}_{n,k} = \frac{1 + \sum_{i=1}^n X_{i,k}\xi_{i,k}}{1 + N_{n,k}}$, $k = 1, \ldots, K$. Here, adding 1 to both the numerator and the denominator is only to avoid the nonsense case of $0/0$. We consider the adaptive design with adding rules $\mathbf{D}_n = \mathbf{D}(\widehat{\boldsymbol{\Theta}}_{n-1}, \boldsymbol{\xi}_n)$ and generating matrices $\mathbf{H}_n = \mathbf{H}(\widehat{\boldsymbol{\Theta}}_{n-1}) = \mathsf{E}[\mathbf{D}_n|\mathcal{F}_{n-1}]$.

The above SEU model has generating matrices $\mathbf{H}_n$ that are not homogeneous. Bai and Hu [4, 5] established the consistency and asymptotic normality of centered urn compositions $\mathbf{Y}_n - \mathsf{E}\mathbf{Y}_n$ for a wide spectrum of urn models with generating matrices $\mathbf{H}_n$ satisfying the condition that

$$\sum_{j=1}^{\infty} \frac{\|\mathbf{H}_j - \mathbf{H}\|}{j} < \infty. \tag{2.1}$$

Here $\|\cdot\|$ denotes the norm of a matrix or a vector. This condition is in general true for practical applications and it is true for the SEU model.

Now, let us turn to the asymptotic property of $\mathbf{N}_n$ that is more useful but intricate. Janson [17] gave a general result of the limiting distribution of $\mathbf{N}_n$ for homogeneous urn models (cf. his Theorem 3.28). The basic technique of Janson [17] is to embed the urn process in a continuous-time branching process. The embedding technique was developed by Athreya and Karlin [3]. Unfortunately, this method does not work for nonhomogeneous urns. The very recent work of Bai and Hu [5] established the asymptotic normality of $\mathbf{N}_n - \mathsf{E}\mathbf{N}_n$ for nonhomogeneous urn models when the generating matrices $\mathbf{H}_n$ satisfy condition (2.1) as well as the following condition:

$$\|\mathbf{H}_n - \mathsf{E}\mathbf{H}_n\| = o(n^{-1/2}) \qquad \text{in } L_2. \tag{2.2}$$

It is beyond doubt that the asymptotic distributions of $\mathbf{Y}_n$ and $\mathbf{N}_n$ with the centralization factors $\mathsf{E}\mathbf{Y}_n$ and $\mathsf{E}\mathbf{N}_n$ being substituted by $n\mathbf{v}$ are more realistic, because the values $\mathsf{E}\mathbf{Y}_n$ and $\mathsf{E}\mathbf{N}_n$ are both unknown and computing them is far too complex, if not infeasible. To provide a possible solution, Bai and Hu [5] suggested that $\mathsf{E}\mathbf{Y}_n$ and $\mathsf{E}\mathbf{N}_n$ can be replaced by $n\mathbf{v}$ if condition (2.1) is revised to a more stringent condition:

$$\sum_{j=1}^{\infty} \frac{\|\mathbf{H}_j - \mathbf{H}\|}{\sqrt{j}} < \infty. \tag{2.3}$$

The above condition is too strict for usual applications. For example, condition (2.3) is not true for the SEU model. On the other hand, condition (2.1) is satisfied with the SEU model because

$$\mathbf{H}(\widehat{\mathbf{\Theta}}_{n-1}) - \mathbf{H}(\mathbf{\Theta}) \approx \sum_{k=1}^{n} \frac{\partial \mathbf{H}(\mathbf{x})}{\partial x_k}\bigg|_{\mathbf{\Theta}} (\hat{\theta}_{n-1,k} - \theta_k) = o(n^{-1/4}) \qquad \text{a.s.}$$

as indicated in Bai and Hu [5]. In fact, conditions (2.2) and (2.3) are not satisfied with the SEU model because the fastest almost sure convergence rate of $\|\mathbf{H}_n - \mathbf{H}\|$ is $O(\sqrt{(\log\log n)/n})$ according to the law of the iterated logarithm and the fastest $L_2$ convergence rate is $O(\sqrt{n})$ according to the central limit theorem. So, the general theorems of Bai and Hu [5] do not apply here.

2.2. *Targeting pre-specified limiting allocation proportions.* As explained in the Introduction, it is frequent that desired limiting allocation proportions are given and our SEU model has the ability to target these proportions. First, let us state the well-known results related to the convergence of $\mathbf{N}_n$ and the generating matrix $\mathbf{H}$ (see, e.g., [4, 5, 17]). Under suitable conditions the allocation proportion converges to a limit, that is,

$$\text{(2.4)} \qquad \frac{\mathbf{N}_n}{n} \to \mathbf{v} \qquad \text{a.s.,}$$

where $\mathbf{v} = (v_1, \ldots, v_K)$ is the left eigenvector of the (limiting) generating matrix $\mathbf{H}$ with respect to its largest eigenvalue and satisfying $v_1 + \cdots + v_K = 1$, that is,

$$\text{(2.5)} \qquad \mathbf{v}\mathbf{H} = \lambda_{\max}\mathbf{v}.$$

To target a pre-specified allocation proportion, one needs to define the SEU process in such a way to force the sample allocation proportion to the desired one as a limit. To illustrate, we will take the two-treatment case with binary responses as an example.

For a binary response (success and failure) clinical trial with two treatments (1 and 2), Rosenberger et al. [26] proposed the following allocation proportion (called the optimal allocation proportion) by considering simultaneously the statistical testing power for comparing the two treatments and the failure rate of patients undergoing the two treatments,

$$\frac{N_{n,1}}{n} \doteq \frac{\sqrt{p_1}}{\sqrt{p_1} + \sqrt{p_2}}, \qquad \frac{N_{n,2}}{n} \doteq \frac{\sqrt{p_2}}{\sqrt{p_1} + \sqrt{p_2}}.$$

Other allocation proportions and a general procedure to find an allocation proportion according to an optimization criterion had been discussed by Jennison and Turnbull [18]. When the SEU model is used to target the above proportion, we need to define an urn model such that $\mathbf{v} = (\frac{\sqrt{p_1}}{\sqrt{p_1}+\sqrt{p_2}}, \frac{\sqrt{p_2}}{\sqrt{p_1}+\sqrt{p_2}})$ is the eigenvector of its generating matrix $\mathbf{H}$ with respect to the largest eigenvalue. Assume that the expectation of total number of balls added at each stage is the same, say $\gamma$. Then $\mathbf{H}$ can be written in the following form:

$$\mathbf{H} = \begin{pmatrix} \gamma - \alpha & \alpha \\ \beta & \gamma - \beta \end{pmatrix}. \tag{2.6}$$

According to (2.5), we have

$$\frac{\alpha}{\beta} = \frac{v_2}{v_1} = \frac{\sqrt{p_2}}{\sqrt{p_1}}.$$

One simple solution of the above equation is $\alpha = \sqrt{p_2}$ and $\beta = \sqrt{p_1}$. Another obvious solution is $\alpha = v_2$ and $\beta = v_1$. In either case, the generating matrix is a function of $p_1$ and $p_2$. If $p_1$ and $p_2$ are known, to implement the allocation proportion, we can define an SEU model with addition rules $\mathbf{D}_i \equiv \mathbf{H}$. Further, when $p_1$ and $p_2$ are unknown, their estimates will be employed as instructed in the previous subsection. Examples are provided in Section 4.

**3. Asymptotic properties.** In this section we will study the asymptotic properties for our SEU model. Major theorems will be given with proofs stated in the last section. Owing to the use of estimates of the unknown parameters, one will discover that the asymptotic variance is different from that given in [5] and [17] and the mathematical techniques are much more complex.

Recall that $\mathbf{Y}_n$ and $\mathbf{N}_n$ are defined by (1.1) and (1.2) and $\mathbf{H} = \mathbf{H}(\boldsymbol{\Theta})$, $\mathbf{D}_n = \mathbf{D}(\widehat{\boldsymbol{\Theta}}_{n-1}, \boldsymbol{\xi}_n)$. To study the asymptotic properties of $\mathbf{Y}_n$ and $\mathbf{N}_n$, the following assumptions are required.

ASSUMPTION 3.1. We assume that $\mathbf{H}\mathbf{1}' = \gamma \mathbf{1}'$, where $\gamma > 0$ and $\mathbf{1} = (1, \ldots, 1)$. Suppose $\mathbf{H}$ has the following Jordan decomposition:

$$\mathbf{T}^{-1}\mathbf{H}\mathbf{T} = \gamma \operatorname{diag}[1, \mathbf{J}_2, \ldots, \mathbf{J}_s],$$

where $\mathbf{J}_s$ is a $\nu_t \times \nu_t$ matrix, given by

$$\mathbf{J}_t = \begin{pmatrix} \lambda_t & 1 & 0 & \ldots & 0 \\ 0 & \lambda_t & 1 & \ldots & 0 \\ 0 & 0 & \lambda_t & \ldots & 0 \\ \vdots & \vdots & \vdots & \ddots & \vdots \\ 0 & 0 & 0 & \ldots & \lambda_t \end{pmatrix}.$$

Without loss of generality, we assume that $\gamma = 1$. Otherwise, we can consider $\mathbf{Y}_n/\gamma$ and $\mathbf{D}_n/\gamma$, instead.

Let $\mathbf{v}$ be the left eigenvector of $\mathbf{H}$ corresponding to its maximal eigenvalue 1 with $\mathbf{v}\mathbf{1}' = 1'$. We may select the matrix $\mathbf{T}$ so that its first column is $\mathbf{1}'$ and the first row of $\mathbf{T}^{-1}$ is $\mathbf{v}$. Let $\lambda = \max\{\text{Re}(\lambda_2), \ldots, \text{Re}(\lambda_s)\}$ and $\nu = \max_j\{\nu_j : \text{Re}(\lambda_j) = \lambda\}$.

ASSUMPTION 3.2. For the addition rules $\mathbf{D}_n = \mathbf{D}(\widehat{\boldsymbol{\Theta}}_{n-1}, \boldsymbol{\xi}_n)$ and the responses $\boldsymbol{\xi}_n$, suppose $D_{ij}(\mathbf{x}, \mathbf{y}) \geq 0$ is continuous at point $(\boldsymbol{\Theta}, \mathbf{y}_0)$ for each $\mathbf{y}_0$ and there is an $r > 2$ such that $\mathsf{E}\|\boldsymbol{\xi}_1\|^r < \infty$ and

$$\sup_n \mathsf{E}[\|\mathbf{D}_n\|^r | \mathcal{F}_{n-1}] < \infty \quad \text{a.s.} \tag{3.1}$$

ASSUMPTION 3.3. For the function $\mathbf{H}(\mathbf{x})$,

$$\mathbf{H}(\mathbf{x}) - \mathbf{H} = \mathbf{H}(\mathbf{x}) - \mathbf{H}(\boldsymbol{\Theta}) = O(\|\mathbf{x} - \boldsymbol{\Theta}\|) \quad \text{as } \mathbf{x} \to \boldsymbol{\Theta}.$$

ASSUMPTION 3.4. For the function $\mathbf{H}(\mathbf{x})$, there is a $\delta > 0$ such that

$$\mathbf{H}(\mathbf{x}) - \mathbf{H} = \sum_{k=1}^K \frac{\partial \mathbf{H}(\mathbf{x})}{\partial x_k}\bigg|_{\mathbf{x}=\boldsymbol{\Theta}} (x_k - \theta_k) + O(\|\mathbf{x} - \boldsymbol{\Theta}\|^{1+\delta}) \quad \text{as } \mathbf{x} \to \boldsymbol{\Theta}.$$

REMARK 3.1. Condition (3.1) is satisfied if we assume that $\mathsf{E}[\|\mathbf{D}(\mathbf{x}, \boldsymbol{\xi}_1)\|^r]$ is a continuous function of $\mathbf{x}$.

Using the above notation, we shall establish the following results.

THEOREM 3.1. *For the SEU, suppose that Assumption 3.1 is satisfied. If $\lambda < 1$, and Assumption 3.2 is satisfied for some $r > 1$, then*

$$\frac{\mathbf{Y}_n}{n} \to \mathbf{v} \quad \text{and} \quad \frac{\mathbf{N}_n}{n} \to \mathbf{v} \quad \text{a.s.}$$

THEOREM 3.2. *For the SEU, if Assumptions 3.1–3.3 are satisfied, and $v_k > 0$, $k = 1, \ldots, K$, $\lambda < 1$, then for any $\kappa > (1/2) \vee \lambda$,*

$$n^{-\kappa}(\mathbf{Y}_n - n\mathbf{v}) \to \mathbf{0} \quad \text{and} \quad n^{-\kappa}(\mathbf{N}_n - n\mathbf{v}) \to \mathbf{0} \quad \text{a.s.}$$

THEOREM 3.3. *For the SEU, if Assumptions 3.1, 3.2 and 3.4 are satisfied, and $v_k > 0$, $k = 1, \ldots, K$, $\lambda < 1/2$, then*

$$n^{1/2}(\mathbf{Y}_n/n - \mathbf{v}) \xrightarrow{\mathcal{D}} N(\mathbf{0}, \mathbf{\Lambda}^{\dagger}) \quad \text{and} \quad n^{1/2}(\mathbf{N}_n/n - \mathbf{v}) \xrightarrow{\mathcal{D}} N(\mathbf{0}, \mathbf{\Lambda}^{\sharp}),$$

*where $N(\mathbf{0}, \mathbf{\Lambda}^{\dagger})$ and $N(\mathbf{0}, \mathbf{\Lambda}^{\sharp})$ are the normal distributions in $\mathcal{R}^K$ with mean zero and $K \times K$ variance–covariance matrices $\mathbf{\Lambda}^{\dagger}$ and $\mathbf{\Lambda}^{\sharp}$, which will be specified by (5.24) and (5.26), respectively.*

Let $\overline{\mathbf{H}} = \mathbf{H} - \mathbf{1}'\mathbf{v}$ and $\mathbf{D} = \mathbf{D}(\boldsymbol{\Theta}, \boldsymbol{\xi}_1)$. Define $\boldsymbol{\Sigma}_1 = \mathrm{diag}(\mathbf{v}) - \mathbf{v}'\mathbf{v}$, $\boldsymbol{\Sigma}_2 = \mathsf{E}[(\mathbf{D}-\mathbf{H})' \mathrm{diag}(\mathbf{v})(\mathbf{D}-\mathbf{H})]$, $\boldsymbol{\Sigma}_3 = \mathrm{diag}(v_1\sigma_1^2, \ldots, v_K\sigma_K^2)$, $\boldsymbol{\Sigma}_{23} = \mathsf{E}[(\mathbf{D}-\mathbf{H})' \times \mathrm{diag}(\mathbf{v})\,\mathrm{diag}(\boldsymbol{\xi}_1 - \boldsymbol{\Theta})]$, where $\sigma_k^2 = \mathrm{Var}(\xi_{i,k})$, $k=1,2,\ldots,K$. Covariance–variance matrices $\boldsymbol{\Lambda}^\dagger$ and $\boldsymbol{\Lambda}^\sharp$ are functions of $\overline{\mathbf{H}}$, $\boldsymbol{\Sigma}_1$, $\boldsymbol{\Sigma}_2$, $\boldsymbol{\Sigma}_3$ and $\boldsymbol{\Sigma}_{23}$.

REMARK 3.2. If $\lambda = 1/2$, the asymptotic normalities also hold, but with different normalizing, given by $n^{1/2}\log^{\nu-1/2} n$, and different $\boldsymbol{\Lambda}^\dagger$ and $\boldsymbol{\Lambda}^\sharp$.

From the proof of Theorem 3.3, one also obtains the following corollary.

COROLLARY 3.1. *For SEU, if Assumptions* 3.1 *and* 3.2 *are satisfied, and* $v_k > 0$, $k=1,\ldots,K$, $\lambda < 1$, *then*

$$n^{1/2}(\widehat{\boldsymbol{\Theta}}_n - \boldsymbol{\Theta}) \xrightarrow{\mathcal{D}} N(\mathbf{0}, \mathrm{diag}(\sigma_1^2/v_1, \ldots, \sigma_K^2/v_K)),$$

*where* $\sigma_k^2 = \mathrm{Var}(\xi_{1,k})$, $k=1,\ldots,K$.

The next corollary provides results on a special case of SEU.

COROLLARY 3.2. *Let* $\boldsymbol{\rho}(\mathbf{x}) = (\rho_K(\mathbf{x}), \ldots, \rho_K(\mathbf{x}))$ *be a vector function being twice differentiable at* $\boldsymbol{\Theta}$ *and* $\rho_k(\mathbf{x}) > 0$, $h=1,\ldots,K$. *Consider an adaptive design with addition rules* $\mathbf{D}_n = \mathbf{1}'\boldsymbol{\rho}(\widehat{\boldsymbol{\Theta}}_{n-1})$, *that is,* $\rho_k(\widehat{\boldsymbol{\Theta}}_{n-1})$ *balls of type* $k$ *are added to the urn at stage* $n$, $k=1,\ldots,K$, *where* $\widehat{\boldsymbol{\Theta}}_{n-1}$ *is the sample estimate of* $\boldsymbol{\Theta}$ *which is defined as in SEU. Then*

$$\frac{\mathbf{Y}_n}{n} \to \boldsymbol{\rho}(\boldsymbol{\Theta}) \quad a.s., \qquad \frac{\mathbf{N}_n}{n} \to \mathbf{v} \quad a.s.$$

*and*

$$\sqrt{n}\left(\frac{\mathbf{Y}_n}{n} - \boldsymbol{\rho}(\boldsymbol{\Theta})\right) \xrightarrow{D} N(\mathbf{0}, 2\boldsymbol{\Sigma}_\rho),$$

$$\sqrt{n}\left(\frac{\mathbf{N}_n}{n} - \mathbf{v}\right) \xrightarrow{D} N(\mathbf{0}, \mathrm{diag}(\mathbf{v}) - \mathbf{v}'\mathbf{v} + 6\boldsymbol{\Sigma}_v),$$

*where* $\gamma = \sum_k \rho_k(\boldsymbol{\Theta})$, $\mathbf{v} = \boldsymbol{\rho}(\boldsymbol{\Theta})/\gamma$,

$$\boldsymbol{\Sigma}_\rho = \left(\frac{\partial \boldsymbol{\rho}(\boldsymbol{\Theta})}{\partial \boldsymbol{\Theta}}\right)' \mathbf{I}^{-1}(\boldsymbol{\Theta}) \frac{\partial \boldsymbol{\rho}(\boldsymbol{\Theta})}{\partial \boldsymbol{\Theta}},$$

$$\boldsymbol{\Sigma}_v = \left(\frac{\partial \mathbf{v}(\boldsymbol{\Theta})}{\partial \boldsymbol{\Theta}}\right)' \mathbf{I}^{-1}(\boldsymbol{\Theta}) \frac{\partial \mathbf{v}(\boldsymbol{\Theta})}{\partial \boldsymbol{\Theta}},$$

$$\mathbf{I}^{-1}(\boldsymbol{\Theta}) = \mathrm{diag}\left(\frac{\sigma_1^2}{v_1}, \ldots, \frac{\sigma_K^2}{v_K}\right)$$

*and* $\mathbf{v}(\mathbf{x}) = \frac{\boldsymbol{\rho}(\mathbf{x})}{\sum_k \rho_k(\mathbf{x})}$. *Here* $\frac{\partial \boldsymbol{\rho}(\boldsymbol{\Theta})}{\partial \boldsymbol{\Theta}}$ *denotes the matrix* $\{\frac{\partial \rho_j(\mathbf{x})}{\partial x_i}|_{\mathbf{x}=\boldsymbol{\Theta}}; i=1,\ldots,K, j=1,\ldots,K\}$.

**4. Applications and conclusions.** In this section we give some motivating examples in order to show the applicability of our SEU model. The first one is related to allocation schemes which aim at assigning more patients to the better treatment. The remaining two are connected to designs with focus on desired targets.

EXAMPLE 1 ([6]). Consider a $K$-treatment clinical trial with binary responses. Let $p_k$ be the success probability of a patient on treatment $k$, and let $q_k = 1 - p_k$ be the failure probability, $k = 1, \ldots, K$. Bai, Hu and Shen [6] (BHS Adaptive Design) proposed to add particles proportional to the estimated success rates for other $K - 1$ treatments, when the response is a failure on treatment $k$, that is, at the $n$th stage, a success on treatment $k$ generates a particle of type $k$, and a failure on treatment $k$ generates $\frac{\hat{p}_{n-1,j}}{(M_{n-1} - \hat{p}_{n-1,k})}$ particles of type $j$ for all $j \neq k$. Here $M_{n-1} = \sum_{j=1}^{K} \hat{p}_{n-1,j}$, $\hat{p}_{n-1,j} = \frac{S_{n-1,j}+1}{N_{n-1,j}+1}$, and $S_{n-1,j}$ denotes the number of successes of treatment $j$ in all the $N_{n-1,j}$ trials of the previous $n-1$ stages, $j = 1, \ldots, K$. They established the strong consistency of $\mathbf{Y}_n$ and $\mathbf{N}_n$ of this model. For this design, $H(\mathbf{x}) = (h_{k,j}(\mathbf{x}), k, j = 1, \ldots, K)$, where $h_{k,k}(\mathbf{x}) = p_k$ and $h_{k,j}(\mathbf{x}) = q_k \frac{x_j}{\sum_{i \neq k} x_i}$ for $k \neq j$. Using Theorem 3.3, we can get the asymptotic normalities.

EXAMPLE 2. We consider a two-treatment clinical trial with binary responses. Suppose that we want to target the optimal allocation proportion proposed by Rosenberger et al. [26]. In (2.6), we choose $\alpha = \sqrt{p_2}$, $\beta = \sqrt{p_1}$ and $\gamma = \sqrt{p_1} + \sqrt{p_2}$. The design is defined as follows. At the $n$th stage, no matter what the response of the $n$th patient is, we add $\sqrt{\hat{p}_{n-1,1}}$ particles of type 1 and $\sqrt{\hat{p}_{n-1,2}}$ particles of type 2 to the urn. If $0 < p_1, p_2 < 1$, then

$$\mathbf{D}(\mathbf{x}) = \mathbf{H}(\mathbf{x}) = \begin{pmatrix} \sqrt{x_1} & \sqrt{x_2} \\ \sqrt{x_1} & \sqrt{x_2} \end{pmatrix} \to \mathbf{D} = \mathbf{H} = \begin{pmatrix} \sqrt{p_1} & \sqrt{p_2} \\ \sqrt{p_1} & \sqrt{p_2} \end{pmatrix}.$$

It can be verified that

$$\mathbf{\Sigma}_\rho = \mathrm{diag}\left(\frac{q_1\sqrt{p_1}}{4(\sqrt{p_1}+\sqrt{p_2})}, \frac{q_2\sqrt{p_2}}{4(\sqrt{p_1}+\sqrt{p_2})}\right)$$

and

$$\mathbf{\Sigma}_v = \frac{1}{4(\sqrt{p_1}+\sqrt{p_2})^3}\left(\frac{p_2 q_1}{\sqrt{p_1}} + \frac{p_1 q_2}{\sqrt{p_2}}\right)\begin{pmatrix} 1 & -1 \\ -1 & 1 \end{pmatrix}.$$

By Corollary 3.2,

$$\left(\frac{Y_{n,1}}{n}, \frac{Y_{n,2}}{n}\right) \to (\sqrt{p_1}, \sqrt{p_2}),$$

$$\left(\frac{N_{n,1}}{n}, \frac{N_{n,2}}{n}\right) \to \left(\frac{\sqrt{p_1}}{\sqrt{p_1}+\sqrt{p_2}}, \frac{\sqrt{p_2}}{\sqrt{p_1}+\sqrt{p_2}}\right) \quad \text{a.s.},$$

$$n^{1/2}\left(\frac{Y_{n,1}}{n} - \sqrt{p_1}, \frac{Y_{n,2}}{n} - \sqrt{p_2}\right) \xrightarrow{\mathcal{D}} N(\mathbf{0}, \mathbf{\Lambda}^\dagger)$$

and

$$n^{1/2}\left(\frac{N_{n,1}}{n} - \frac{\sqrt{p_1}}{\sqrt{p_1} + \sqrt{p_2}}\right) \xrightarrow{\mathcal{D}} N(0, \sigma_\sharp^2),$$

where

$$\mathbf{\Lambda}^\dagger = \operatorname{diag}\left(\frac{q_1\sqrt{p_1}}{2(\sqrt{p_1} + \sqrt{p_2})}, \frac{q_2\sqrt{p_2}}{2(\sqrt{p_1} + \sqrt{p_2})}\right)$$

and

$$\sigma_\sharp^2 = \frac{\sqrt{p_1 p_2}}{(\sqrt{p_1} + \sqrt{p_2})^2} + \frac{3}{2(\sqrt{p_1} + \sqrt{p_2})^3}\left(\frac{p_2 q_1}{\sqrt{p_1}} + \frac{p_1 q_2}{\sqrt{p_2}}\right).$$

EXAMPLE 3. We consider a two-treatment clinical trial with binary responses. Suppose that we want to target the allocation proportion ($v_1 = \frac{q_2}{q_1+q_2}, v_2 = \frac{q_1}{q_1+q_2}$) of the randomized play-the-winner (RPW) rule proposed by Wei and Durham [31]. In (2.6) we choose $\alpha = v_2$, $\beta = v_1$ and $\gamma = 1$. Then $\mathbf{H} = (1,1)'(v_1, v_2)$. The design is defined as follows. At the $n$th stage, regardless of what the response of the $n$th patient is, we add $\frac{\hat{q}_{n-1,2}}{\hat{q}_{n-1,1}+\hat{q}_{n-1,2}}$ particles of type 1 and $\frac{\hat{q}_{n-1,1}}{\hat{q}_{n-1,1}+\hat{q}_{n-1,2}}$ particles of type 2 to the urn, where $\hat{q}_{n-1,k} = 1 - \hat{p}_{n-1,k}$, $k = 1, 2$. If $0 < p_1, p_2 < 1$, then

$$\mathbf{D}(\mathbf{x}) = \begin{pmatrix} \frac{1-x_2}{(1-x_1)+(1-x_2)} & \frac{1-x_1}{(1-x_1)+(1-x_2)} \\ \frac{1-x_2}{(1-x_1)+(1-x_2)} & \frac{1-x_1}{(1-x_1)+(1-x_2)} \end{pmatrix}$$

$$\to \mathbf{D} = \begin{pmatrix} \frac{q_2}{q_1+q_2} & \frac{q_1}{q_1+q_2} \\ \frac{q_2}{q_1+q_2} & \frac{q_1}{q_1+q_2} \end{pmatrix},$$

and $\mathbf{H}(\mathbf{x}) = \mathbf{D}(\mathbf{x})$, $\mathbf{H} = \mathbf{D}$. It can be verified that $\mathbf{\Sigma}_\rho = \mathbf{\Sigma}_v = q_1 q_2 (p_1 + p_2)/(q_1+q_2)^3 (1,-1)'(-1,1)$. By Corollary 3.2,

$$\left(\frac{Y_{n,1}}{n}, \frac{Y_{n,2}}{n}\right) \to \left(\frac{q_2}{q_1+q_2}, \frac{q_1}{q_1+q_2}\right),$$

$$\left(\frac{N_{n,1}}{n}, \frac{N_{n,2}}{n}\right) \to \left(\frac{q_2}{q_1+q_2}, \frac{q_1}{q_1+q_2}\right) \quad \text{a.s.,}$$

$$n^{1/2}\left(\frac{Y_{n,1}}{n} - \frac{q_2}{q_1+q_2}, \frac{Y_{n,2}}{n} - \frac{q_1}{q_1+q_2}\right) \xrightarrow{\mathcal{D}} N(0, \sigma_\dagger^2)(1,-1)$$

and

$$n^{1/2}\left(\frac{N_{n,1}}{n} - \frac{q_2}{q_1+q_2}\right) \xrightarrow{\mathcal{D}} N(0, \sigma_\sharp^2),$$

where

$$\sigma_\dagger^2 = \frac{2q_1q_2(p_1+p_2)}{(q_1+q_2)^3} \quad \text{and} \quad \sigma_\sharp^2 = \frac{q_1q_2[2+5(p_1+p_2)]}{(q_1+q_2)^3}.$$

The asymptotic variances are much smaller than those given by the RPW rule when $q_1 + q_2$ is near or less than $1/2$, since when $q_1 + q_2 > 1/2$, the asymptotic variances of $Y_{n,1}/n$ and $N_{n,1}/n$ in the RPW rule are

$$\frac{q_1q_2}{[2(q_1+q_2)-1](q_1+q_2)^2} \quad \text{and} \quad \frac{q_1q_2[1+2(p_1+p_2)]}{[2(q_1+q_2)-1](q_1+q_2)^2},$$

respectively (cf. [28]). When $q_1 + q_2 < 1/2$, the limiting distributions of $\mathbf{Y}_n$ and $\mathbf{N}_n$ in the RPW rule are unknown, and the convergence rates of the variances are much slower. But, for the adaptive design given here, the asymptotic normality can be evaluated according to our formula for any $0 < q_1, q_2 < 1$.

Examples 2 and 3 demonstrate the ability of the SEU model to target a given allocation proportion. In general, we can use Corollary 3.2 to target any desired allocation proportion. This is a useful property of our proposed model.

The examples in this section provide illustrations of how the proposed techniques are translated into applications in the formulation of valuable treatment allocation schemes. We can also evaluate these schemes by reviewing the respective asymptotic properties of $\mathbf{Y}_n$ and $\mathbf{N}_n$.

The main contributions of this paper are in proposing the family of sequential estimation-adjusted urn models and showing the two important properties of the SEU models: (i) they can be used to target any given allocation proportions, and (ii) they have certain desired asymptotic properties under some widely satisfied conditions. Based on the asymptotic results, we are able to derive new designs which have smaller variabilities than some traditional and popular urn models (see Example 3).

In this paper we considered the application of the SEU models in clinical trials. The idea of using constantly updated estimators in SEU models is novel and the SEU models can be employed in other areas. In fact, there are abundant urn model applications where data are collected sequentially and the consecutive estimations of unknown parameters provide better solutions to the formulated problem.

For simplicity, the examples considered in this section are all binary responses. However, the SEU models can be applied to different types of responses. From the definition of the SEU process, the response $\boldsymbol{\xi}_n$ could be discrete or continuous. For example, we can apply the SEU model to the case studied in Section 8 of [12]. In that case, the response is normally distributed and the target proportion is the Neyman allocation. Generally, the assumptions in Section 3 are satisfied if the third moments of $\boldsymbol{\xi}$ and $\mathbf{D}$ exist and the function $\mathbf{H}$ is differentiable. Therefore, the SEU model has a wide spectrum of applications.

**5. Proofs.** First, we state several relevant lemmas. Recall $\overline{\mathbf{H}} = \mathbf{H} - \mathbf{1}'\mathbf{v}$. Then

$$(5.1) \qquad \mathbf{T}^{-1}\overline{\mathbf{H}}\mathbf{T} = \mathrm{diag}[0, \mathbf{J}_2, \ldots, \mathbf{J}_s].$$

LEMMA 5.1. *Let $\overline{\mathbf{B}}_{n,n} = \mathbf{I}$ and $\overline{\mathbf{B}}_{n,i} = \prod_{j=i+1}^{n}(\mathbf{I} + j^{-1}\overline{\mathbf{H}})$. If two sequences of matrices $\mathbf{Q}_n$ and $\mathbf{P}_n$ satisfy $\mathbf{Q}_0 = \mathbf{P}_0 = \mathbf{0}$ and*

$$\mathbf{Q}_n = \mathbf{P}_n + \sum_{k=0}^{n-1} \frac{\mathbf{Q}_k}{k+1}\overline{\mathbf{H}},$$

*that is,*

$$\mathbf{Q}_n = \Delta\mathbf{P}_n + \mathbf{Q}_{n-1}(\mathbf{I} + n^{-1}\overline{\mathbf{H}}),$$

*where $\Delta\mathbf{P}_n = \mathbf{P}_n - \mathbf{P}_{n-1}$ is the difference of $\mathbf{P}_n$, then*

$$(5.2) \qquad \mathbf{Q}_n = \sum_{m=1}^{n} \Delta\mathbf{P}_m \overline{\mathbf{B}}_{n,m} = \mathbf{P}_n + \sum_{m=1}^{n-1} \mathbf{P}_m \frac{\overline{\mathbf{H}}}{m+1} \overline{\mathbf{B}}_{n,m+1}.$$

*Also,*

$$(5.3) \quad \|\overline{\mathbf{B}}_{n,m}\| \leq C(n/m)^\lambda \log^{\nu-1}(n/m) \qquad \textit{for all } m = 1, \ldots, n, n \geq 1,$$

*where $\log x = \ln(x \vee e)$ hereafter.*

PROOF. See Lemma A.1 of [16]. □

LEMMA 5.2. *For the SEU, if Assumption 3.2 is satisfied with some $r > 1$, then $N_{n,k} \to \infty$, $k = 1, \ldots, K$, as $n \to \infty$. Furthermore, $\widehat{\boldsymbol{\Theta}}_n \to \boldsymbol{\Theta}$ a.s. whenever $\mathsf{E}\|\boldsymbol{\xi}_1\| < \infty$, and $\hat{\theta}_{n,k} - \theta_k = O(\sqrt{\frac{\log\log N_{n,k}}{N_{n,k}}})$ a.s. whenever $\mathsf{E}\|\boldsymbol{\xi}_1\|^2 < \infty$ as $n \to \infty$, $k = 1, \ldots, K$.*

PROOF. First, notice that $\{\|\mathbf{D}_n\| - \mathsf{E}[\|\mathbf{D}_n\| | \mathcal{F}_{n-1}]\}$ is a sequence of martingale differences with

$$\sup_n \mathsf{E}[\|\mathbf{D}_n\| | \mathcal{F}_{n-1}] \leq \sup_n (\mathsf{E}[\|\mathbf{D}_n\|^r | \mathcal{F}_{n-1}])^{1/r} < \infty \qquad \text{a.s.}$$

and

$$\sum_{n=1}^{\infty} \frac{\mathsf{E}[|\|\mathbf{D}_n\| - \mathsf{E}[\|\mathbf{D}_n\||\mathcal{F}_{n-1}]|^r|\mathcal{F}_{n-1}]}{n^r} \leq \sum_{n=1}^{\infty} \frac{C}{n^r} < \infty \qquad \text{a.s.}$$

by condition (3.1). According to the law of large numbers of martingales (cf. [14]),

$$\sum_{m=1}^{n} \|\mathbf{D}_m\| = \sum_{m=1}^{n} \mathsf{E}[\|\mathbf{D}_m\||\mathcal{F}_{m-1}] + o(n) = O(n) \qquad \text{a.s.}$$

So, by (1.1), $\|\mathbf{Y}_n\| \leq \|\mathbf{Y}_0\| + \sum_{m=1}^{n} \|\mathbf{D}_m\| = O(n)$ a.s. It follows that $\sum_{n=1}^{\infty}(\sum_j Y_{n,j})^{-1} = \infty$ a.s. So, for each $k = 1, \ldots, K$,

$$\sum_{n=1}^{\infty} \mathsf{P}(X_{n,k} = 1|\mathcal{F}_{n-1}) = \sum_{n=1}^{\infty} \frac{Y_{n-1,k}}{\sum_j Y_{n,j}} \geq \sum_{n=1}^{\infty} \frac{Y_{0,k}}{\sum_j Y_{n,j}} = \infty \qquad \text{a.s.},$$

which, together with the generalized Borel–Cantelli lemmas, implies that $\mathsf{P}(X_{n,k} = 1, \text{i.o.}) = 1$. Then $N_{n,k} \to \infty$ a.s., $k = 1, \ldots, K$. By Lemma A.4 of [16] again, the lemma is proved. □

Now, let $|\mathbf{Y}_n| = \mathbf{Y}_n\mathbf{1}'$ be the total number of balls in the urn at the $n$th stage, and let $\mathbf{M}_n = \sum_{k=1}^{n} \Delta\mathbf{M}_k$ and $\mathbf{m}_n = \sum_{k=1}^{n} \Delta\mathbf{m}_k$, where $\Delta\mathbf{M}_k = \mathbf{X}_k\mathbf{D}_k - \mathsf{E}[\mathbf{X}_k\mathbf{D}_k|\mathcal{F}_{k-1}]$ and $\Delta\mathbf{m}_k = \mathbf{X}_k - \mathsf{E}[\mathbf{X}_k|\mathcal{F}_{k-1}]$. Then $\{(\mathbf{m}_n, \mathbf{M}_n), \mathcal{F}_n; n \geq 1\}$ is a martingale sequence. By (1.1),

$$(5.4) \qquad \mathbf{Y}_n = \mathbf{Y}_0 + \mathbf{M}_n + \sum_{m=0}^{n-1} \frac{\mathbf{Y}_m}{|\mathbf{Y}_m|}\mathbf{H}_{m+1}.$$

Note that $\mathbf{v}\mathbf{1}'\mathbf{v} = \mathbf{v}$ and $\mathbf{v}\overline{\mathbf{H}} = \mathbf{0}$. By (5.4), it follows that

$$\begin{aligned}
\mathbf{Y}_n - n\mathbf{v} &= \mathbf{M}_n + \sum_{m=0}^{n-1}\left(\frac{\mathbf{Y}_m}{|\mathbf{Y}_m|} - \mathbf{v}\right)\mathbf{H}_{m+1} + \sum_{m=1}^{n}\mathbf{v}(\mathbf{H}_m - \mathbf{H}) + \mathbf{Y}_0 \\
&= \mathbf{M}_n + \sum_{m=0}^{n-1}\left(\frac{\mathbf{Y}_m}{|\mathbf{Y}_m|} - \mathbf{v}\right)\mathbf{H} + \sum_{m=1}^{n}\mathbf{v}(\mathbf{H}_m - \mathbf{H}) \\
&\quad + \sum_{m=0}^{n-1}\left(\frac{\mathbf{Y}_m}{|\mathbf{Y}_m|} - \mathbf{v}\right)(\mathbf{H}_{m+1} - \mathbf{H}) + \mathbf{Y}_0 \\
(5.5) \\
&= \mathbf{M}_n + \sum_{m=0}^{n-1}\left(\frac{\mathbf{Y}_m}{|\mathbf{Y}_m|} - \mathbf{v}\right)\overline{\mathbf{H}} + \sum_{m=1}^{n}\mathbf{v}(\mathbf{H}_m - \mathbf{H}) \\
&\quad + \sum_{m=0}^{n-1}\left(\frac{\mathbf{Y}_m}{|\mathbf{Y}_m|} - \mathbf{v}\right)(\mathbf{H}_{m+1} - \mathbf{H}) + \mathbf{Y}_0
\end{aligned}$$

$$= \mathbf{M}_n + \sum_{m=0}^{n-1} \frac{\mathbf{Y}_m - m\mathbf{v}}{m+1}\overline{\mathbf{H}} + \mathbf{R}_n,$$

where

(5.6)
$$\mathbf{R}_n = \sum_{m=1}^{n} \mathbf{v}(\mathbf{H}_m - \mathbf{H}) + \sum_{m=0}^{n-1}\left(\frac{\mathbf{Y}_m}{|\mathbf{Y}_m|} - \mathbf{v}\right)\left(1 - \frac{|\mathbf{Y}_m|}{m+1}\right)\overline{\mathbf{H}}$$
$$+ \sum_{m=0}^{n-1}\left(\frac{\mathbf{Y}_m}{|\mathbf{Y}_m|} - \mathbf{v}\right)(\mathbf{H}_{m+1} - \mathbf{H}) + \mathbf{Y}_0.$$

Also by (1.2),

$$\mathbf{N}_n - n\mathbf{v} = \mathbf{m}_n + \sum_{m=0}^{n-1} \frac{\mathbf{Y}_m}{|\mathbf{Y}_m|} - n\mathbf{v}$$

$$= \mathbf{m}_n + \sum_{m=0}^{n-1}\left(\frac{\mathbf{Y}_m}{|\mathbf{Y}_m|} - \mathbf{v}\right)(\mathbf{I} - \mathbf{1}'\mathbf{v})$$

$$= \mathbf{m}_n + \sum_{m=0}^{n-1}\left(\frac{\mathbf{Y}_m}{m+1} - \mathbf{v}\right)(\mathbf{I} - \mathbf{1}'\mathbf{v})$$

(5.7)
$$+ \sum_{m=0}^{n-1} \frac{\mathbf{Y}_m}{|\mathbf{Y}_m|}\left(1 - \frac{|\mathbf{Y}_m|}{m+1}\right)(\mathbf{I} - \mathbf{1}'\mathbf{v})$$

$$= \mathbf{m}_n + \sum_{m=0}^{n-1} \frac{\mathbf{Y}_m - m\mathbf{v}}{m+1}(\mathbf{I} - \mathbf{1}'\mathbf{v})$$

$$+ \sum_{m=0}^{n-1}\left(\frac{\mathbf{Y}_m}{|\mathbf{Y}_m|} - \mathbf{v}\right)\left(1 - \frac{|\mathbf{Y}_m|}{m+1}\right)(\mathbf{I} - \mathbf{1}'\mathbf{v}).$$

PROOF OF THEOREM 3.1. First, note that $r > 1$, $\|\Delta \mathbf{m}_n\| \leq 2$, and according to (3.1),

(5.8) $\quad \mathsf{E}[\|\Delta \mathbf{M}_n\|^r | \mathcal{F}_{n-1}] \leq 2^r \sup_n \mathsf{E}[\|\mathbf{D}_n\|^r | \mathcal{F}_{n-1}] := 2^r \eta_r < \infty.$

By the strong law of large numbers of martingales, we have

(5.9) $\quad\quad\quad\quad \mathbf{M}_n = o(n) \quad \text{and} \quad \mathbf{m}_n = o(n) \quad\quad \text{a.s.}$

On the other hand, by Lemma 5.2 and the continuity of $\mathbf{H}(\mathbf{x})$, we have

(5.10) $\quad\quad\quad\quad\quad\quad \mathbf{H}_n - \mathbf{H} = o(1) \quad\quad \text{a.s.}$

Note that $\overline{\mathbf{H}}\mathbf{1}' = \mathbf{0}'$, and $\mathbf{Y}_m/|\mathbf{Y}_m| - \mathbf{v}$ is bounded. From (5.6), it follows that

$$\mathbf{R}_n\mathbf{1}' = \sum_{m=1}^{n} \mathbf{v}(\mathbf{H}_m - \mathbf{H})\mathbf{1}' + \sum_{m=0}^{n-1}\left(\frac{\mathbf{Y}_m}{|\mathbf{Y}_m|} - \mathbf{v}\right)(\mathbf{H}_{m+1} - \mathbf{H})\mathbf{1}' + \mathbf{Y}_0\mathbf{1}' \quad (5.11)$$
$$= o(n) \quad \text{a.s.,}$$

which, together with (5.5) and (5.9), implies that

$$|\mathbf{Y}_n| - n = (\mathbf{Y}_n - n\mathbf{v})\mathbf{1}' = \mathbf{M}_n\mathbf{1}' + \mathbf{R}_n\mathbf{1}' = o(n) \quad \text{a.s.,}$$

that is,

$$(5.12) \qquad \frac{|\mathbf{Y}_n|}{n} \to 1 \quad \text{a.s.}$$

Now, by combining (5.6), (5.10) and (5.12) it follows that

$$\mathbf{R}_n = o(n) \quad \text{a.s.,}$$

which, together with (5.5) and (5.9), yields

$$\mathbf{Y}_n - n\mathbf{v} = \sum_{m=0}^{n-1} \frac{\mathbf{Y}_m - m\mathbf{v}}{m+1}\overline{\mathbf{H}} + o(n) \quad \text{a.s.}$$

And then by Lemma 5.1, we conclude that

$$\mathbf{Y}_n - n\mathbf{v} = o(n) + \sum_{m=1}^{n-1} \frac{o(m)}{m+1}(n/m)^\lambda \log^{\nu-1}(n/m) = o(n) \quad \text{a.s.,}$$

that is,

$$\frac{\mathbf{Y}_n}{n} \to \mathbf{v} \quad \text{a.s.}$$

Finally, from (5.7) it follows that

$$\mathbf{N}_n - n\mathbf{v} = o(n) + \sum_{m=0}^{n-1} \frac{o(m)}{m+1} + \sum_{m=0}^{n-1}\left(\frac{\mathbf{Y}_m}{|\mathbf{Y}_m|} - \mathbf{v}\right)o(1) = o(n) \quad \text{a.s.}$$

The proof of Theorem 3.1 is complete. □

To prove Theorem 3.2, we need the following lemma.

LEMMA 5.3. *Under the assumptions in Theorem 3.2, we have*

$$\widehat{\mathbf{\Theta}}_n - \mathbf{\Theta} = O\left(\sqrt{\frac{\log\log n}{n}}\right) \quad \text{and} \quad \mathbf{H}_n - \mathbf{H} = O\left(\sqrt{\frac{\log\log n}{n}}\right) \quad a.s.$$

PROOF. Note that by Theorem 3.1,
$$\frac{n}{N_{n,k}} \to \frac{1}{v_k} \qquad \text{a.s., } k=1,\ldots,K.$$
So, by Lemma 5.2,
$$\hat{\theta}_{n,k} - \theta_k = O\left(\sqrt{\frac{\log\log n}{n}}\right) \qquad \text{a.s., } k=1,\ldots,K.$$
The proof is then complete by noting Assumption 3.3. □

PROOF OF THEOREM 3.2. First, note that $\|\Delta \mathbf{m}_n\| \leq 2$, (5.8) and $r>2$. By the law of the iterated logarithm of martingales, we have

(5.13) $\mathbf{M}_n = O(\sqrt{n\log\log n})$ and $\mathbf{m}_n = O(\sqrt{n\log\log n})$ a.s.

By Lemma 5.3 and similarly to (5.11) we have
$$\mathbf{R}_n \mathbf{1}' = \sum_{m=1}^{n} O\left(\sqrt{\frac{\log\log m}{m}}\right) + \sum_{m=1}^{n-1}\left(\frac{\mathbf{Y}_m}{|\mathbf{Y}_m|} - \mathbf{v}\right) O\left(\sqrt{\frac{\log\log m}{m}}\right) \qquad \text{a.s.}$$
$$= O(\sqrt{n\log\log n}),$$
which, together with (5.5), implies that
$$|\mathbf{Y}_n| - n = \mathbf{M}_n \mathbf{1}' + \mathbf{R}_n \mathbf{1}' = O(\sqrt{n\log\log n}) \qquad \text{a.s.,}$$
that is,

(5.14) $$\frac{|\mathbf{Y}_n|}{n} - 1 = O\left(\sqrt{\frac{\log\log n}{n}}\right) \qquad \text{a.s.}$$

And then from (5.6), it follows that
$$\mathbf{R}_n = \sum_{m=1}^{n} O\left(\sqrt{\frac{\log\log m}{m}}\right) + \sum_{m=1}^{n-1}\left(\frac{\mathbf{Y}_m}{|\mathbf{Y}_m|} - \mathbf{v}\right) O\left(\sqrt{\frac{\log\log m}{m}}\right) \qquad \text{a.s.}$$
$$= O(\sqrt{n\log\log n}).$$
With (5.5) and (5.13), we conclude that
$$\mathbf{Y}_n - n\mathbf{v} = \sum_{m=0}^{n-1} \frac{\mathbf{Y}_m - m\mathbf{v}}{m+1} \overline{\mathbf{H}} + O(\sqrt{n\log\log n}) \qquad \text{a.s.,}$$
which, together with Lemma 5.1, implies that
$$\mathbf{Y}_n - n\mathbf{v} = \sum_{m=1}^{n} \frac{O(\sqrt{m\log\log m})}{m+1}(n/m)^\lambda \log^{\nu-1}(n/m) \qquad \text{a.s.}$$

(5.15)
$$= \begin{cases} O(\sqrt{n\log\log n}), & \text{if } \lambda < 1/2, \\ O(\sqrt{n\log\log n}\,\log^\nu n), & \text{if } \lambda = 1/2, \\ O(n^\lambda \log^{\nu-1} n), & \text{if } \lambda > 1/2. \end{cases}$$

Finally, by (5.7) and (5.13)–(5.15) we conclude that

$$\mathbf{N}_n - n\mathbf{v} = O(\sqrt{n\log\log n}) + \sum_{m=0}^{n-1} \frac{O(\|\mathbf{Y}_m - m\mathbf{v}\|)}{m+1}$$

$$+ \sum_{m=1}^{n-1} O\left(\sqrt{\frac{\log\log m}{m}}\right) \quad \text{a.s.}$$

$$= \begin{cases} O(\sqrt{n\log\log n}), & \text{if } \lambda < 1/2, \\ O(\sqrt{n\log\log n}\log^\nu n), & \text{if } \lambda = 1/2, \\ O(n^\lambda \log^{\nu-1} n), & \text{if } \lambda > 1/2. \end{cases}$$

Theorem 3.2 is proved. □

PROOF OF THEOREM 3.3. Without loss of generality, we assume that $\delta$ is small enough such that $\delta < 1/2 - \lambda$. Define

$$(5.16) \quad \mathbf{f}_k = \mathbf{v}\frac{\partial \mathbf{H}(\mathbf{x})}{\partial x_k}\bigg|_{\mathbf{x}=\boldsymbol{\Theta}}, \quad \mathbf{F} = \begin{pmatrix} \mathbf{f}_1/v_1 \\ \vdots \\ \mathbf{f}_K/v_K \end{pmatrix} = (\mathbf{f}_1'/v_1, \ldots, \mathbf{f}_K'/v_K)'$$

and $\mathbf{Q}_n = (Q_{n,1}, \ldots, Q_{n,K}) = \sum_{m=1}^n \Delta \mathbf{Q}_m$, where $\Delta \mathbf{Q}_m = \mathbf{X}_m \operatorname{diag}(\xi_m - \boldsymbol{\Theta})$. Then $\mathbf{Q}_n$ is a martingale since $\xi_m$ is independent of $\mathbf{X}_m$ and the $\sigma$-field $\mathcal{F}_{m-1}$. By Lemma 5.3 and Assumption 3.4,

$$\mathbf{v}(\mathbf{H}_{n+1} - \mathbf{H}) = (\widehat{\boldsymbol{\Theta}}_n - \boldsymbol{\Theta})(\mathbf{f}_1', \ldots, \mathbf{f}_K')' + O(\|\widehat{\boldsymbol{\Theta}}_n - \boldsymbol{\Theta}\|^{1+\delta})$$
$$= (\widehat{\boldsymbol{\Theta}}_n - \boldsymbol{\Theta})(\mathbf{f}_1', \ldots, \mathbf{f}_K')' + o(n^{-1/2-\delta/3}) \quad \text{a.s.}$$

Also, by Theorem 3.2, $\mathbf{Q}_n = O(\sqrt{n\log\log n})$ a.s. and $\hat{\theta}_{n,k} = \frac{1+\sum_{m=1}^n X_{m,k}\xi_{m,k}}{1+N_{n,k}}$, $k = 1, \ldots, K$, we have

$$\hat{\theta}_{n,k} - \theta_k = \frac{n}{N_{n,k}+1}\frac{1}{n}(Q_{n,k} + 1 - \theta_k)$$

$$= \frac{1}{v_k}\frac{1}{n}Q_{n,k} + \left(\frac{n}{N_{n,k}+1} - \frac{1}{v_k}\right)\frac{1}{n}Q_{n,k} + \frac{1-\theta_k}{N_{n,k}+1}$$

$$= \frac{1}{v_k}\frac{1}{n}Q_{n,k} + o(n^{\kappa-1})O(n^{-1/2}\sqrt{\log\log n}) + O\left(\frac{1}{n}\right)$$

$$= \frac{1}{v_k}\frac{1}{n}Q_{n,k} + o(n^{-1/2-\delta}) \quad \text{a.s.}$$

It follows that

$$(5.17) \quad \mathbf{v}(\mathbf{H}_{n+1} - \mathbf{H}) = \frac{1}{n}\mathbf{Q}_n\mathbf{F} + o(n^{-1/2-\delta/3}) \quad \text{a.s.}$$

Now, from (5.6), (5.17), Theorem 3.2 and Lemma 5.3 it follows that

$$\mathbf{R}_n = \sum_{m=1}^{n} \mathbf{v}(\mathbf{H}_m - \mathbf{H}) + \sum_{m=1}^{n-1} o((m^{-1+\kappa})^2)$$

$$+ \sum_{m=1}^{n-1} o(m^{-1+\kappa}) O\left(\sqrt{\frac{\log\log m}{m}}\right)$$

$$= \sum_{m=1}^{n} \frac{1}{m} \mathbf{Q}_m \mathbf{F} + o(n^{1/2-\delta/3}) \quad \text{a.s.}$$

So, by (5.5),

(5.18)
$$\mathbf{Y}_n - n\mathbf{v} = \mathbf{M}_n + \sum_{m=1}^{n} \frac{1}{m} \mathbf{Q}_m \mathbf{F} + \sum_{m=0}^{n-1} \frac{\mathbf{Y}_m - m\mathbf{v}}{m+1} \overline{\mathbf{H}}$$
$$+ o(n^{-1/2-\delta/3}) \quad \text{a.s.}$$

Let $\mathbf{U}_0 = \mathbf{0}$ and

$$\mathbf{U}_n = \mathbf{M}_n + \sum_{m=1}^{n} \frac{1}{m} \mathbf{Q}_m \mathbf{F} + \sum_{m=0}^{n-1} \frac{\mathbf{U}_m}{m+1} \overline{\mathbf{H}}.$$

Then

$$\mathbf{Y}_n - n\mathbf{v} - \mathbf{U}_n = \sum_{m=0}^{n-1} \frac{\mathbf{Y}_m - m\mathbf{v} - \mathbf{U}_m}{m+1} \overline{\mathbf{H}} + o(n^{1/2-\delta/3}) \quad \text{a.s.}$$

By Lemma 5.1, it follows that

(5.19)
$$\mathbf{Y}_n - n\mathbf{v} = \mathbf{U}_n + \sum_{m=1}^{n} o(m^{-1/2-\delta/3})\left(\frac{n}{m}\right)^{\lambda} \log^{\nu-1}\left(\frac{n}{m}\right)$$
$$= \mathbf{U}_n + o((n^{1/2-\delta/4}) \vee (n^{\lambda} \log^{\nu-1} n)) \quad \text{a.s.}$$

On the other hand, by (5.2),

$$\mathbf{U}_n = \sum_{m=1}^{n} \left(\Delta \mathbf{M}_m + \frac{1}{m} \mathbf{Q}_m \mathbf{F}\right) \overline{\mathbf{B}}_{n,m}$$
$$= \sum_{m=1}^{n} \Delta \mathbf{M}_m \overline{\mathbf{B}}_{n,m} + \sum_{m=1}^{n} \Delta \mathbf{Q}_m \mathbf{F} \sum_{j=m}^{n} \left(\frac{1}{j} \overline{\mathbf{B}}_{n,j}\right).$$

Also, by (5.7) and Theorem 3.2,

$$\mathbf{N}_n - n\mathbf{v} = \mathbf{m}_n + \sum_{m=0}^{n-1} \frac{\mathbf{Y}_m - m\mathbf{v}}{m+1}(\mathbf{I} - \mathbf{1}'\mathbf{v}) + \sum_{m=1}^{n-1} o((m^{-1+\kappa})^2)$$

$$= \mathbf{m}_n + \sum_{m=0}^{n-1} \frac{\mathbf{U}_m}{m+1}(\mathbf{I} - \mathbf{1}'\mathbf{v})$$
$$+ o((n^{1/2-\delta/4}) \vee (n^\lambda \log^{\nu-1} n)) + o(n^{1/2-\delta})$$

(5.20)
$$= \mathbf{m}_n + \sum_{m=1}^{n-1} \Delta\mathbf{M}_m \left(\sum_{i=m}^{n-1} \frac{1}{i+1}\overline{\mathbf{B}}_{i,m}\right)(\mathbf{I} - \mathbf{1}'\mathbf{v})$$
$$+ \sum_{m=1}^{n-1} \Delta\mathbf{Q}_m \mathbf{F}\left(\sum_{i=m}^{n-1}\sum_{j=m}^{i} \frac{1}{(i+1)j}\overline{\mathbf{B}}_{i,j}\right)(\mathbf{I} - \mathbf{1}'\mathbf{v})$$
$$+ o((n^{1/2-\delta/4}) \vee (n^\lambda \log^{\nu-1} n))$$
$$:= \mathbf{V}_n + o((n^{1/2-\delta/4}) \vee (n^\lambda \log^{\nu-1} n)) \quad \text{a.s.}$$

So, to show the asymptotic normalities of $\mathbf{Y}_n$ and $\mathbf{N}_n$, it suffices to show the asymptotic normalities of $\mathbf{U}_n$ and $\mathbf{V}_n$. Write $\mathbf{B}_{n,m}^{(1)} = \sum_{j=m}^{n} \frac{\overline{\mathbf{B}}_{n,j}}{j}$, $\mathbf{B}_{n,m}^{(2)} = \sum_{i=m}^{n-1} \frac{\overline{\mathbf{B}}_{i,m}}{i+1}$ and $\mathbf{B}_{n,m}^{(3)} = \sum_{i=m}^{n-1}\sum_{j=m}^{i} \frac{\overline{\mathbf{B}}_{i,j}}{(i+1)j}$, where $\sum_{j=k+1}^{k}(\cdot) = 0$ and so $\mathbf{B}_{n,n}^{(i)} = \mathbf{0}$, $i = 1,2,3$. Then

(5.21)
$$\mathbf{U}_n = \sum_{m=1}^{n}(\Delta\mathbf{M}_m \overline{\mathbf{B}}_{n,m} + \Delta\mathbf{Q}_m \mathbf{F}\mathbf{B}_{n,m}^{(1)}),$$
$$\mathbf{V}_n = \sum_{m=1}^{n}(\Delta\mathbf{m}_m + \Delta\mathbf{M}_m \mathbf{B}_{n,m}^{(2)}(\mathbf{I} - \mathbf{1}'\mathbf{v}) + \Delta\mathbf{Q}_m \mathbf{F}\mathbf{B}_{n,m}^{(3)}(\mathbf{I} - \mathbf{1}'\mathbf{v})).$$

Note that $\mathbf{U}_n$ and $\mathbf{V}_n$ are sums of martingale differences. We will use the central limit theorem for martingale (cf. [14]) to prove our result. According to (5.8),

$$\frac{1}{n^{r/2}} \sum_{m=1}^{n} \{\mathsf{E}[\|\Delta\mathbf{M}_m \overline{\mathbf{B}}_{n,m} + \Delta\mathbf{Q}_m \mathbf{F}\mathbf{B}_{n,m}^{(1)}\|^r | \mathcal{F}_{m-1}]$$
$$+ \mathsf{E}[\|\Delta\mathbf{m}_m + \Delta\mathbf{M}_m \mathbf{B}_{n,m}^{(2)}(\mathbf{I} - \mathbf{1}'\mathbf{v})$$
$$+ \Delta\mathbf{Q}_m \mathbf{F}\mathbf{B}_{n,m}^{(3)}(\mathbf{I} - \mathbf{1}'\mathbf{v})\|^r | \mathcal{F}_{m-1}]\}$$

(5.22)
$$\leq \frac{C}{n^{r/2}} \sum_{m=1}^{n} \left(\left(\frac{n}{m}\right)^\lambda \log^{\nu-1} \frac{n}{m}\right)^r$$
$$\times \mathsf{E}[\|\Delta\mathbf{m}_m\|^r + \|\Delta\mathbf{M}_m\|^r + \|\Delta\mathbf{Q}_m\|^r | \mathcal{F}_{m-1}]$$
$$\leq \frac{C}{n^{r/2}} \sum_{m=1}^{n} \left(\left(\frac{n}{m}\right)^\lambda \log^{\nu-1} \frac{n}{m}\right)^r \{1 + 2^r \eta_r + \mathsf{E}\|\boldsymbol{\xi}_m\|^r\}$$

$$\leq \frac{Cn}{n^{r/2}} \int_1^n y^{\lambda r - 2} \log^{(\nu-1)r} y \, dy \to 0 \qquad \text{a.s. as } n \to \infty,$$

that is, the Lindberg condition is satisfied. Finally, it is enough to calculate the asymptotic (conditioned) variance–covariance matrices of $\mathbf{U}_n$ and $\mathbf{V}_n$. Recall $\mathbf{\Sigma}_1 = \text{diag}(\mathbf{v}) - \mathbf{v}'\mathbf{v}$, $\mathbf{\Sigma}_2 = \mathsf{E}[(\mathbf{D} - \mathbf{H})' \text{diag}(\mathbf{v})(\mathbf{D} - \mathbf{H})]$, $\mathbf{\Sigma}_3 = \text{diag}(v_1 \sigma_1^2, \ldots, v_K \sigma_K^2)$, $\mathbf{\Sigma}_{23} = \mathsf{E}[(\mathbf{D} - \mathbf{H})' \text{diag}(\mathbf{v}) \text{diag}(\boldsymbol{\xi}_1 - \boldsymbol{\Theta})]$, where $\mathbf{D} = \mathbf{D}(\boldsymbol{\Theta}, \boldsymbol{\xi}_1)$ and $\sigma_k^2 = \text{Var}(\xi_{i,k})$, $k = 1, 2, \ldots, K$. Notice (5.8) and the continuity of the functions $\mathbf{D}(\cdot, \boldsymbol{\xi}_n)$ and $\mathbf{H}(\cdot)$. Then by Theorem 3.1, Assumption 3.2 and Lemma 5.2, we have

$$\text{Var}[\Delta \mathbf{M}_n | \mathcal{F}_{n-1}]$$
$$= \mathsf{E}\left[\mathbf{D}_n' \text{diag}\left(\frac{\mathbf{Y}_{n-1}}{|\mathbf{Y}_{n-1}|}\right) \mathbf{D}_n \Big| \mathcal{F}_{n-1}\right] - \mathbf{H}_n' \frac{\mathbf{Y}_{n-1}'}{|\mathbf{Y}_{n-1}|} \frac{\mathbf{Y}_{n-1}}{|\mathbf{Y}_{n-1}|} \mathbf{H}_n$$
$$= \mathsf{E}\left[(\mathbf{D}_n - \mathbf{H}_n)' \text{diag}\left(\frac{\mathbf{Y}_{n-1}}{|\mathbf{Y}_{n-1}|}\right) (\mathbf{D}_n - \mathbf{H}_n) \Big| \mathcal{F}_{n-1}\right]$$
$$+ \mathbf{H}_n' \left[\text{diag}\left(\frac{\mathbf{Y}_{n-1}}{|\mathbf{Y}_{n-1}|}\right) - \frac{\mathbf{Y}_{n-1}'}{|\mathbf{Y}_{n-1}|} \frac{\mathbf{Y}_{n-1}}{|\mathbf{Y}_{n-1}|}\right] \mathbf{H}_n$$
$$\to \mathbf{\Sigma}_2 + \mathbf{H}' \mathbf{\Sigma}_1 \mathbf{H} \qquad \text{a.s.,}$$

$$\text{Var}[\Delta \mathbf{m}_n | \mathcal{F}_{n-1}]$$
$$= \text{diag}\left(\frac{\mathbf{Y}_{n-1}}{|\mathbf{Y}_{n-1}|}\right) - \frac{\mathbf{Y}_{n-1}'}{|\mathbf{Y}_{n-1}|} \frac{\mathbf{Y}_{n-1}}{|\mathbf{Y}_{n-1}|} \to \mathbf{\Sigma}_1 \qquad \text{a.s.,}$$

$$\text{Var}[\Delta \mathbf{Q}_n | \mathcal{F}_{n-1}]$$
$$= \mathsf{E}[\text{diag}(X_{n,1}(\xi_{n,1} - \theta_1)^2, \ldots, X_{n,K}(\xi_{n,K} - \theta_K)^2) | \mathcal{F}_{n-1}]$$
$$= \text{diag}\left(\frac{\sigma_1^2 Y_{n-1,1}}{|\mathbf{Y}_{n-1}|}, \ldots, \frac{\sigma_K^2 Y_{n-1,K}}{|\mathbf{Y}_{n-1}|}\right) \to \mathbf{\Sigma}_3 \qquad \text{a.s.,}$$

$$\text{Cov}[(\Delta \mathbf{m}_n, \Delta \mathbf{M}_n) | \mathcal{F}_{n-1}]$$
$$= \left[\text{diag}\left(\frac{\mathbf{Y}_{n-1}}{|\mathbf{Y}_{n-1}|}\right) - \frac{\mathbf{Y}_{n-1}'}{|\mathbf{Y}_{n-1}|} \frac{\mathbf{Y}_{n-1}}{|\mathbf{Y}_{n-1}|}\right] \mathbf{H}_n$$
$$\to \mathbf{\Sigma}_1 \mathbf{H} \qquad \text{a.s.,}$$

$$\text{Cov}[(\Delta \mathbf{m}_n, \Delta \mathbf{Q}_n) | \mathcal{F}_{n-1}]$$
$$= \mathsf{E}[(\mathbf{X}_n - \mathsf{E}[\mathbf{X}_n | \mathcal{F}_{n-1}])' \mathbf{X}_n \text{diag}(\boldsymbol{\xi}_n - \boldsymbol{\Theta}) | \mathcal{F}_{n-1}]$$
$$= \mathbf{0}$$

and

$$\text{Cov}[(\Delta \mathbf{M}_n, \Delta \mathbf{Q}_n) | \mathcal{F}_{n-1}]$$

$$= \mathsf{E}\left[\mathbf{D}'_n \operatorname{diag}\left(\frac{\mathbf{Y}_{n-1}}{|\mathbf{Y}_{n-1}|}\right) \operatorname{diag}(\boldsymbol{\xi}_n - \boldsymbol{\Theta})|\mathcal{F}_{n-1}\right]$$
$$\to \boldsymbol{\Sigma}_{23} \quad \text{a.s.}$$

So,

$$\operatorname{Var}[(\Delta \mathbf{m}_n, \Delta \mathbf{M}_n - \Delta \mathbf{m}_n \mathbf{H}, \Delta \mathbf{Q}_n)|\mathcal{F}_{n-1}]$$

(5.23)
$$\to \begin{pmatrix} \boldsymbol{\Sigma}_1 & \mathbf{0} & \mathbf{0} \\ \mathbf{0} & \boldsymbol{\Sigma}_2 & \boldsymbol{\Sigma}_{23} \\ \mathbf{0} & \boldsymbol{\Sigma}'_{23} & \boldsymbol{\Sigma}_3 \end{pmatrix} \quad \text{a.s.}$$

By noting (5.21),

$$\mathbf{H}\overline{\mathbf{B}}_{n,m} = \overline{\mathbf{B}}_{n,m}\mathbf{H} \quad \text{and} \quad \mathbf{H}(\mathbf{I} - \mathbf{1}'\mathbf{v}) = \overline{\mathbf{H}},$$

we can write

$$\mathbf{U}_n = \mathbf{U}_n^{(1)}\mathbf{H} + \mathbf{U}_n^{(2)} + \mathbf{U}_n^{(3)}$$

and

$$\mathbf{V}_n = \mathbf{V}_n^{(1)} + \mathbf{V}_n^{(2)}(\mathbf{I} - \mathbf{1}'\mathbf{v}) + \mathbf{V}_n^{(3)}(\mathbf{I} - \mathbf{1}'\mathbf{v}),$$

where

$$\mathbf{U}_n^{(1)} = \sum_{m=1}^{n} \Delta \mathbf{m}_m \overline{\mathbf{B}}_{n,m},$$

$$\mathbf{U}_n^{(2)} = \sum_{m=1}^{n} (\Delta \mathbf{M}_m - \Delta \mathbf{m}_n \mathbf{H})\overline{\mathbf{B}}_{n,m},$$

$$\mathbf{U}_n^{(3)} = \sum_{m=1}^{n} \Delta \mathbf{Q}_m \mathbf{F} \mathbf{B}_{n,m}^{(1)},$$

$$\mathbf{V}_n^{(1)} = \sum_{m=1}^{n} \Delta \mathbf{m}_m + \sum_{m=1}^{n-1} \Delta \mathbf{m}_m \mathbf{H} \mathbf{B}_{n,m}^{(2)}(\mathbf{I} - \mathbf{1}'\mathbf{v})$$

$$= \Delta \mathbf{m}_n + \sum_{m=1}^{n-1} \Delta \mathbf{m}_m (\mathbf{B}_{n,m}^{(2)} \overline{\mathbf{H}} + \mathbf{I})$$

$$= \Delta \mathbf{m}_n + \sum_{m=1}^{n-1} \Delta \mathbf{m}_m \overline{\mathbf{B}}_{n,m} = \mathbf{U}_n^{(1)},$$

$$\mathbf{V}_n^{(2)} = \sum_{m=1}^{n} (\Delta \mathbf{M}_m - \Delta \mathbf{m}_m \mathbf{H}) \mathbf{B}_{n,m}^{(2)},$$

$$\mathbf{V}_n^{(3)} = \sum_{m=1}^{n} \Delta \mathbf{Q}_m \mathbf{F} \mathbf{B}_{n,m}^{(3)}.$$

Then by (5.23),

$$\operatorname{Var}(\mathbf{U}_n^{(1)}|) = \sum_{m=1}^n \mathsf{E}[(\Delta \mathbf{m}_m \overline{\mathbf{B}}_{n,m})' \Delta \mathbf{m}_m \overline{\mathbf{B}}_{n,m} | \mathcal{F}_{m-1}]$$

$$= \sum_{m=1}^n \overline{\mathbf{B}}'_{n,m}(\mathbf{\Sigma}_1 + o(1))\overline{\mathbf{B}}_{n,m}$$

$$= \int_1^n \left(\frac{n}{x}\right)^{\overline{\mathbf{H}}'} \mathbf{\Sigma}_1 \left(\frac{n}{x}\right)^{\overline{\mathbf{H}}} dx + o(1) \sum_{m=1}^n \left(\frac{n}{m}\right)^{2\lambda} \log^{2\nu-2}\left(\frac{n}{m}\right)$$

$$= n \int_{1/n}^1 \left(\frac{1}{x}\right)^{\overline{\mathbf{H}}'} \mathbf{\Sigma}_1 \left(\frac{1}{x}\right)^{\overline{\mathbf{H}}} dx + o(n)$$

$$= n \int_0^1 \left(\frac{1}{x}\right)^{\overline{\mathbf{H}}'} \mathbf{\Sigma}_1 \left(\frac{1}{x}\right)^{\overline{\mathbf{H}}} dx + o(n) := n \mathbf{\Lambda}_1^\dagger + o(n) \quad \text{a.s.},$$

where $a^{\overline{\mathbf{H}}}$ is defined to be $e^{\overline{\mathbf{H}} \ln a} = \sum_{j=0}^\infty \frac{(\ln a)^j}{j!} \overline{\mathbf{H}}^j$. Here we use $\operatorname{Var}(\cdot|)$ to denote the sum of conditional variance–covariance matrices of related martingale differences, and $\operatorname{Cov}\{\cdot,\cdot|\}$ is defined similarly. Also,

$$\operatorname{Var}(\mathbf{U}_n^{(2)}|) = \sum_{m=1}^n \overline{\mathbf{B}}'_{n,m}(\mathbf{\Sigma}_2 + o(1))\overline{\mathbf{B}}_{n,m}$$

$$= n \int_0^1 \left(\frac{1}{x}\right)^{\overline{\mathbf{H}}'} \mathbf{\Sigma}_2 \left(\frac{1}{x}\right)^{\overline{\mathbf{H}}} dx + o(n) := n \mathbf{\Lambda}_2^\dagger + o(n) \quad \text{a.s.},$$

$$\operatorname{Var}(\mathbf{U}_n^{(3)}|) = \sum_{m=1}^n (\mathbf{B}_{n,m}^{(1)})' \mathbf{F}'(\mathbf{\Sigma}_3 + o(1))\mathbf{F} \mathbf{B}_{n,m}^{(1)}$$

$$= \int_1^n dx \left[\int_x^n \frac{1}{y}\left(\frac{n}{y}\right)^{\overline{\mathbf{H}}} dy\right]' \mathbf{F}'\mathbf{\Sigma}_3 \mathbf{F} \left[\int_x^n \frac{1}{y}\left(\frac{n}{y}\right)^{\overline{\mathbf{H}}} dy\right] + o(n)$$

$$= n \int_0^1 dx \left[\int_x^1 \frac{1}{y}\left(\frac{1}{y}\right)^{\overline{\mathbf{H}}} dy\right]' \mathbf{F}'\mathbf{\Sigma}_3 \mathbf{F} \left[\int_x^1 \frac{1}{y}\left(\frac{1}{y}\right)^{\overline{\mathbf{H}}} dy\right] + o(n)$$

$$:= n \mathbf{\Lambda}_3^\dagger + o(n) \quad \text{a.s.},$$

$$\operatorname{Cov}\{\mathbf{U}_n^{(1)}, \mathbf{U}_n^{(2)}|\} = \sum_{m=1}^n \mathsf{E}[(\Delta \mathbf{m}_m \overline{\mathbf{B}}_{n,m})'(\Delta \mathbf{M}_m - \Delta \mathbf{m}_n \mathbf{H})\overline{\mathbf{B}}_{n,m}|\mathcal{F}_{m-1}]$$

$$= o(n) \quad \text{a.s.},$$

$$\operatorname{Cov}\{\mathbf{U}_n^{(1)}, \mathbf{U}_n^{(3)}|\} = 0 \quad \text{a.s.}$$

and

$$\text{Cov}\{\mathbf{U}_n^{(2)}, \mathbf{U}_n^{(3)}|\} = \sum_{m=1}^n \overline{\mathbf{B}}'_{n,m}(\mathbf{\Sigma}_{23} + o(1))\mathbf{F}\mathbf{B}_{n,m}^{(1)}$$

$$= \int_1^n dx \left(\frac{n}{x}\right)^{\overline{\mathbf{H}}'} \mathbf{\Sigma}_{23}\mathbf{F}\left[\int_x^n \frac{1}{y}\left(\frac{n}{y}\right)^{\overline{\mathbf{H}}} dy\right] + o(n)$$

$$= n\int_0^1 dx \left(\frac{1}{x}\right)^{\overline{\mathbf{H}}'} \mathbf{\Sigma}_{23}\mathbf{F}\left[\int_x^1 \frac{1}{y}\left(\frac{1}{y}\right)^{\overline{\mathbf{H}}} dy\right] + o(n)$$

$$:= n\mathbf{\Lambda}_{23}^\dagger + o(n) \quad \text{a.s.}$$

Note that $\|a^{\overline{\mathbf{H}}}\| \leq Ca^\lambda \log^{\nu-1} a$ for $a \geq 1$. The above integrals are well defined. It follows that

(5.24) $\quad n^{-1}\text{Var}(\mathbf{U}_n|) \to \mathbf{H}'\mathbf{\Lambda}_1^\dagger\mathbf{H} + \mathbf{\Lambda}_2^\dagger + \mathbf{\Lambda}_3^\dagger + \mathbf{\Lambda}_{23}^\dagger + (\mathbf{\Lambda}_{23}^\dagger)' := \mathbf{\Lambda}^\dagger.$

And then by the central limit theorem for martingales (cf. [14])

(5.25) $\quad\quad\quad\quad\quad n^{-1/2}\mathbf{U}_n \xrightarrow{\mathcal{D}} N(\mathbf{0}, \mathbf{\Lambda}^\dagger).$

Similarly,

$$\text{Var}(\mathbf{V}_n^{(1)}|) = \text{Var}(\mathbf{U}_n^{(1)}|) = n\mathbf{\Lambda}_1^\dagger + o(n) \quad \text{a.s.},$$

$$\text{Var}(\mathbf{V}_n^{(2)}|) = \sum_{m=1}^n (\mathbf{B}_{n,m}^{(2)})'(\mathbf{\Sigma}_2 + o(1))\mathbf{B}_{n,m}^{(2)}$$

$$= \int_1^n dx\left[\int_x^n \frac{1}{y}\left(\frac{y}{x}\right)^{\overline{\mathbf{H}}} dy\right]' \mathbf{\Sigma}_2 \left[\int_x^n \frac{1}{y}\left(\frac{y}{x}\right)^{\overline{\mathbf{H}}} dy\right] + o(n)$$

$$= n\int_0^1 dx\left[\int_x^1 \frac{1}{y}\left(\frac{y}{x}\right)^{\overline{\mathbf{H}}} dy\right]' \mathbf{\Sigma}_2 \left[\int_x^1 \frac{1}{y}\left(\frac{y}{x}\right)^{\overline{\mathbf{H}}} dy\right] + o(n)$$

$$:= n\mathbf{\Lambda}_2^\sharp + o(n) \quad \text{a.s.},$$

$$\text{Var}(\mathbf{V}_n^{(3)}|) = \sum_{m=1}^n (\mathbf{B}_{n,m}^{(3)})'\mathbf{F}'(\mathbf{\Sigma}_3 + o(1))\mathbf{F}B_{n,m}^{(3)}$$

$$= \int_1^n dx\left[\int_x^n dy \int_x^y du \frac{1}{yu}\left(\frac{y}{u}\right)^{\overline{\mathbf{H}}}\right]'$$

$$\times \mathbf{F}'\mathbf{\Sigma}_3\mathbf{F}\left[\int_x^n dy \int_x^y du \frac{1}{yu}\left(\frac{y}{u}\right)^{\overline{\mathbf{H}}}\right] + o(n)$$

$$= n\int_0^1 dx\left[\int_x^1 dy \int_x^y du \frac{1}{yu}\left(\frac{y}{u}\right)^{\overline{\mathbf{H}}}\right]'$$

$$\times \mathbf{F}'\boldsymbol{\Sigma}_3 \mathbf{F}\left[\int_x^1 dy \int_x^y du \frac{1}{yu}\left(\frac{y}{u}\right)^{\overline{\mathbf{H}}}\right] + o(n)$$

$$:= n\boldsymbol{\Lambda}_3^\sharp + o(n) \quad \text{a.s.},$$

$$\mathrm{Cov}\{\mathbf{V}_n^{(1)}, \mathbf{V}_n^{(2)}|\} = o(n) \quad \text{a.s.}, \quad \mathrm{Cov}\{\mathbf{V}_n^{(1)}, \mathbf{V}_n^{(3)}|\} = 0 \quad \text{a.s.}$$

and

$$\mathrm{Cov}\{\mathbf{V}_n^{(2)}, \mathbf{V}_n^{(3)}|\}$$

$$= \sum_{m=1}^n (\mathbf{B}_{n,m}^{(2)})'(\boldsymbol{\Sigma}_{23} + o(1))\mathbf{F} B_{n,m}^{(3)}$$

$$= \int_1^n dx \left[\int_x^n \frac{1}{y}\left(\frac{y}{x}\right)^{\overline{\mathbf{H}}} dy\right]' \mathbf{F}'\boldsymbol{\Sigma}_{23}\mathbf{F}\left[\int_x^n dy \int_x^y du \frac{1}{yu}\left(\frac{y}{u}\right)^{\overline{\mathbf{H}}}\right] + o(n)$$

$$= n\int_0^1 dx \left[\int_x^1 \frac{1}{y}\left(\frac{y}{x}\right)^{\overline{\mathbf{H}}} dy\right]' \mathbf{F}'\boldsymbol{\Sigma}_{23}\mathbf{F}\left[\int_x^1 dy \int_x^y du \frac{1}{yu}\left(\frac{y}{u}\right)^{\overline{\mathbf{H}}}\right] + o(n)$$

$$:= n\boldsymbol{\Lambda}_{23}^\sharp + o(n) \quad \text{a.s.}$$

It follows that

$$n^{-1}\mathrm{Var}(\mathbf{V}_n|) \to \boldsymbol{\Lambda}_1^\dagger + (\mathbf{I} - \mathbf{v}'\mathbf{1})\boldsymbol{\Lambda}_2^\sharp(\mathbf{I} - \mathbf{1}'\mathbf{v})$$

(5.26)
$$+ (\mathbf{I} - \mathbf{v}'\mathbf{1})(\boldsymbol{\Lambda}_3^\sharp + \boldsymbol{\Lambda}_{23}^\sharp + (\boldsymbol{\Lambda}_{23}^\sharp)')(\mathbf{I} - \mathbf{1}'\mathbf{v})$$

$$:= \boldsymbol{\Lambda}^\sharp \quad \text{a.s.}$$

And then by the central limit theorem for martingales (cf. [14]),

(5.27) $$n^{-1/2}\mathbf{V}_n \xrightarrow{\mathcal{D}} N(\mathbf{0}, \boldsymbol{\Lambda}^\sharp).$$

Combining (5.19), (5.20), (5.25) and (5.27), we complete the proof of Theorem 3.3. □

PROOF OF COROLLARY 3.1. Note that $\mathbf{Q}_n = O(n^{1/2})$ in $L_2$. By Theorem 3.1, we have

$$\hat{\theta}_{n,k} - \theta_k = \frac{1}{v_k}\frac{1}{n}Q_{n,k} + \left(\frac{n}{N_{n,k}+1} - \frac{1}{v_k}\right)\frac{1}{n}Q_{n,k} + \frac{1-\theta_k}{N_{n,k}+1}$$

$$= \frac{1}{v_k}\frac{1}{n}Q_{n,k} + o(1)O(n^{-1/2}) + O\left(\frac{1}{n}\right)$$

$$= \frac{1}{v_k}\frac{1}{n}Q_{n,k} + o(n^{-1/2}) \quad \text{in probability, } k = 1,\ldots,K.$$

On the other hand,

$$n^{-1/2}\mathbf{Q}_n \xrightarrow{\mathcal{D}} N(\mathbf{0}, \boldsymbol{\Sigma}_3).$$

The result follows.  □

PROOF OF COROLLARY 3.2. Now, $\mathbf{H1}' = \gamma$, we shall consider $\mathbf{Y}_n/\gamma$ instead of $\mathbf{Y}_n$, and $\mathbf{D}(\mathbf{x})/\gamma$ instead of $\mathbf{D}(\mathbf{x})$. Notice that $\mathbf{D}(\mathbf{x})/\gamma \equiv \mathbf{H}(\mathbf{x})/\gamma = \mathbf{1}'\boldsymbol{\rho}(\mathbf{x})/\gamma$, and then $\boldsymbol{\Sigma}_2 = \boldsymbol{\Sigma}_{23} = \mathbf{0}$. It follows that $\boldsymbol{\Lambda}_2^\dagger = \boldsymbol{\Lambda}_{23}^\dagger = \boldsymbol{\Lambda}_2^\sharp = \boldsymbol{\Lambda}_{23}^\sharp = \mathbf{0}$, and then $\boldsymbol{\Lambda}^\dagger = (\mathbf{H}/\gamma)'\boldsymbol{\Lambda}_1^\dagger(\mathbf{H}/\gamma) + \boldsymbol{\Lambda}_3^\dagger$ and $\boldsymbol{\Lambda}^\sharp = \boldsymbol{\Lambda}_1^\dagger + (\mathbf{I} - \mathbf{v}'\mathbf{1})\boldsymbol{\Lambda}_3^\sharp(\mathbf{I} - \mathbf{1}'\mathbf{v})$ by (5.25) and (5.26). Also, $\mathbf{H}/\gamma = \mathbf{1}'\mathbf{v}$ and $\overline{\mathbf{H}} = \mathbf{H}/\gamma - \mathbf{1}'\mathbf{v} = \mathbf{0}$. Thus

$$\boldsymbol{\Lambda}_1^\dagger = \int_0^1 \left(\frac{1}{x}\right)^0 \boldsymbol{\Sigma}_1 \left(\frac{1}{x}\right)^0 dx = \boldsymbol{\Sigma}_1 = \mathrm{diag}(\mathbf{v}) - \mathbf{v}'\mathbf{v},$$

$(\mathbf{H}'/\gamma)\boldsymbol{\Lambda}_1^\dagger(\mathbf{H}/\gamma) = \mathbf{0}$,

$$\boldsymbol{\Lambda}_3^\dagger = \int_0^1 dx \left[\int_x^1 \frac{dy}{y}\right]^2 (\mathbf{F}/\gamma)'\boldsymbol{\Sigma}_3(\mathbf{F}/\gamma) = 2(\mathbf{F}/\gamma)'\boldsymbol{\Sigma}_3(\mathbf{F}/\gamma)$$

and

$$\boldsymbol{\Lambda}_3^\sharp = \int_0^1 dx \left[\int_x^1 \frac{dy}{y} \int_x^y \frac{du}{u}\right]^2 (\mathbf{F}/\gamma)'\boldsymbol{\Sigma}_3(\mathbf{F}/\gamma)$$
$$= \int_0^1 \frac{1}{4}(\ln x)^4 (\mathbf{F}/\gamma)'\boldsymbol{\Sigma}_3(\mathbf{F}/\gamma)\, dx = 6(\mathbf{F}/\gamma)'\boldsymbol{\Sigma}_3(\mathbf{F}/\gamma).$$

Notice that $\mathbf{v}\partial\mathbf{H}(\mathbf{x})/\partial x_k = \partial\boldsymbol{\rho}(\mathbf{x})/\partial x_k$. It is easy to see that

$$(\mathbf{F}/\gamma)'\boldsymbol{\Sigma}_3(\mathbf{F}/\gamma) = \frac{1}{\gamma^2}\left(\frac{\partial\boldsymbol{\rho}(\boldsymbol{\Theta})}{\partial\boldsymbol{\Theta}}\right)' \mathbf{I}^{-1}(\boldsymbol{\Theta})\frac{\partial\boldsymbol{\rho}(\boldsymbol{\Theta})}{\partial\boldsymbol{\Theta}} = \boldsymbol{\Sigma}_\rho/\gamma^2.$$

Hence $\boldsymbol{\Lambda}^\dagger = 2\boldsymbol{\Sigma}_\rho/\gamma^2$ and $\boldsymbol{\Lambda}^\sharp = \boldsymbol{\Lambda}_1^\dagger + 6(\mathbf{I} - \mathbf{v}'\mathbf{1})\boldsymbol{\Sigma}_\rho(\mathbf{I} - \mathbf{1}'\mathbf{v})/\gamma^2$. By Theorems 3.1 and 3.3, we conclude that

$$\frac{\mathbf{Y}_n/\gamma}{n} \to \mathbf{v} \quad \text{a.s.}, \qquad \frac{\mathbf{N}_n}{n} \to \mathbf{v} \quad \text{a.s.}$$

and

$$\sqrt{n}\left(\frac{\mathbf{Y}_n/\gamma}{n} - \mathbf{v}\right) \xrightarrow{D} N(\mathbf{0}, 2\boldsymbol{\Sigma}_\rho/\gamma^2), \qquad \sqrt{n}\left(\frac{\mathbf{N}_n}{n} - \mathbf{v}\right) \xrightarrow{D} N(\mathbf{0}, \boldsymbol{\Lambda}^\sharp).$$

Finally, notice that $\boldsymbol{\rho}(x) = \rho(x)\mathbf{1}'\mathbf{v}(\mathbf{x})$. We have

$$\frac{\partial\boldsymbol{\rho}(\boldsymbol{\Theta})}{\partial\boldsymbol{\Theta}}(\mathbf{I} - \mathbf{1}'\mathbf{v}) = \frac{\partial\rho(\boldsymbol{\Theta})\mathbf{1}'}{\partial\boldsymbol{\Theta}}\mathbf{v}(\mathbf{I} - \mathbf{1}'\mathbf{v}) + \gamma\frac{\partial\mathbf{v}(\boldsymbol{\Theta})}{\partial\boldsymbol{\Theta}}(\mathbf{I} - \mathbf{1}'\mathbf{v})$$
$$= \mathbf{0} + \gamma\frac{\partial\mathbf{v}(\boldsymbol{\Theta})}{\partial\boldsymbol{\Theta}} - \gamma\frac{\partial\mathbf{v}(\boldsymbol{\Theta})\mathbf{1}'}{\partial\boldsymbol{\Theta}}\mathbf{v} = \gamma\frac{\partial\mathbf{v}(\boldsymbol{\Theta})}{\partial\boldsymbol{\Theta}}.$$

The proof is now complete.  □

REMARK 5.1. One can also show that $n^{-1}\operatorname{Cov}\{\mathbf{U}_n, \mathbf{V}_n|\} \to \mathbf{\Lambda}^{\dagger\sharp}$ a.s. for some $\mathbf{\Lambda}^{\dagger\sharp}$. So,

$$n^{1/2}(\mathbf{Y}_n/n - \mathbf{v}, \mathbf{N}_n/n - \mathbf{v}) \overset{\mathcal{D}}{\to} N\left(\mathbf{0}, \begin{pmatrix} \mathbf{\Lambda}^{\dagger} & \mathbf{\Lambda}^{\dagger\sharp} \\ (\mathbf{\Lambda}^{\dagger\sharp})' & \mathbf{\Lambda}^{\sharp} \end{pmatrix}\right).$$

Actually, with a more careful but complex proof similar to that in [16], one can show that the process $n^{-1/2}(\mathbf{Y}_{[nt]} - [nt]\mathbf{v}, \mathbf{N}_{[nt]} - [nt]\mathbf{v})$ is weakly convergent to a $2K$-dimensional Gaussian process $(\mathbf{G}_1(t), \mathbf{G}_2(t))$ in the Skorohod topology, where $\mathbf{G}_1(t)$ is the solution of

$$d\mathbf{G}_1(t) = d\mathbf{W}_1(t)\mathbf{H} + d\mathbf{W}_2(t) + \frac{\mathbf{W}_3(t)}{t}\mathbf{F}\,dt + \frac{\mathbf{G}_1(t)}{t}\overline{\mathbf{H}}\,dt,$$

$$\mathbf{G}_1(0) = \mathbf{0},$$

$(\mathbf{W}_1(t), \mathbf{W}_2(t), \mathbf{W}_3(t))$ is a $3K$-dimensional Brownian motion with variance–covariance matrix defined in (5.23) and

$$\mathbf{G}_2(t) = \mathbf{W}_1(t) + \int_0^t \frac{\mathbf{G}_1(s)}{s}\,ds(\mathbf{I} - \mathbf{1}'\mathbf{v}).$$

Furthermore, we also have

$$n^{-1}\operatorname{Cov}\{\mathbf{Q}_n, \mathbf{U}_n|\} \to (\mathbf{\Sigma}'_{23} + \mathbf{\Sigma}_3)(\mathbf{I} - \overline{\mathbf{H}})^{-1} \qquad \text{a.s.}$$

and

$$n^{-1}\operatorname{Cov}\{\mathbf{Q}_n, \mathbf{V}_n|\} \to (\mathbf{\Sigma}_{23}' + \mathbf{\Sigma}_3)(\mathbf{I} - \overline{\mathbf{H}})^{-1}(\mathbf{I} - \mathbf{1}'\mathbf{v}) \qquad \text{a.s.}$$

This shows that the joint distribution of $n^{1/2}(\widehat{\boldsymbol{\Theta}}_n - \boldsymbol{\Theta})$, $n^{1/2}(\mathbf{Y}_n/n - \mathbf{v})$ and $n^{1/2}(\mathbf{N}_n/n - \mathbf{v})$ converges to a mean-zero $3K$-dimensional Gaussian distribution.

If $\lambda = 1/2$, one also can show that (5.22) holds with $(n\log^{2\nu-1} n)^{1+\delta/2}$ replacing $n^{r/2}$, and that $(n\log^{2\nu-1} n)^{-1}\operatorname{Var}\{(\mathbf{U}_n, \mathbf{V}_n)|\}$ converges to a $2d \times 2d$ matrix. This implies that the asymptotic normalities of $\mathbf{Y}_n$ and $\mathbf{N}_n$ also hold.

REMARK 5.2. $\mathbf{\Lambda}^{\dagger}$ and $\mathbf{\Lambda}^{\sharp}$ can be derived. As an example, we provide details to derive the second part of $\mathbf{\Lambda}^{\sharp}$, that is, $(\mathbf{I} - \mathbf{v}'\mathbf{1})\mathbf{\Lambda}_2^{\sharp}(\mathbf{I} - \mathbf{1}'\mathbf{v})$. Other similar terms can be derived in the same fashion. Write $\mathbf{T} = (\mathbf{1}', \mathbf{T}_2, \ldots, \mathbf{T}_s)$ and $\mathbf{T}^*(\mathbf{I} - \mathbf{v}'\mathbf{1})\mathbf{\Lambda}_2^{\sharp}(\mathbf{I} - \mathbf{1}'\mathbf{v})\mathbf{T} = (\mathbf{\Sigma}_{gh}, g, h = 1, \ldots, s)$. Note that

$\mathbf{T}^*\mathbf{\Lambda}_2^{\sharp}\mathbf{T}$

$$= \int_0^1 dx \left[\int_x^1 \frac{1}{y}\left(\frac{y}{x}\right)^{\operatorname{diag}[0,\mathbf{J}_2^*,\ldots,\mathbf{J}_s^*]} dy\right] \mathbf{T}^*\mathbf{\Sigma}_2\mathbf{T} \left[\int_x^1 \frac{1}{y}\left(\frac{y}{x}\right)^{\operatorname{diag}[0,\mathbf{J}_2,\ldots,\mathbf{J}_s]} dy\right]$$

$$= \int_0^1 dx \left[\int_x^1 \operatorname{diag}\left[1, \frac{1}{y}\left(\frac{y}{x}\right)^{\mathbf{J}_2^*}, \ldots, \frac{1}{y}\left(\frac{y}{x}\right)^{\mathbf{J}_s^*}\right] dy\right] \mathbf{T}^*\mathbf{\Sigma}_2\mathbf{T}$$

$$\times \left[\int_x^1 \operatorname{diag}\left[1, \frac{1}{y}\left(\frac{y}{x}\right)^{\mathbf{J}_2}, \ldots, \frac{1}{y}\left(\frac{y}{x}\right)^{\mathbf{J}_s}\right] dy\right].$$

So,
$$\mathbf{T}^*(\mathbf{I} - \mathbf{v}'\mathbf{1})\mathbf{\Lambda}_2^\sharp(\mathbf{I} - \mathbf{1}'\mathbf{v})\mathbf{T}$$
$$= \mathrm{diag}(0,1,\ldots,1)\mathbf{T}^*\mathbf{\Lambda}_2^\sharp\mathbf{T}\,\mathrm{diag}(0,1,\ldots,1)$$
$$= \int_0^1 dx \left[\int_x^1 \mathrm{diag}\left[0, \frac{1}{y}\left(\frac{y}{x}\right)^{\mathbf{J}_2^*}, \ldots, \frac{1}{y}\left(\frac{y}{x}\right)^{\mathbf{J}_s^*}\right] dy\right]\mathbf{T}^*\mathbf{\Sigma}_2\mathbf{T}$$
$$\times \left[\int_x^1 \mathrm{diag}\left[0, \frac{1}{y}\left(\frac{y}{x}\right)^{\mathbf{J}_2}, \ldots, \frac{1}{y}\left(\frac{y}{x}\right)^{\mathbf{J}_s}\right] dy\right].$$

Also,
$$\left(\frac{y}{x}\right)^{\mathbf{J}_t} = \left(\frac{y}{x}\right)^{\lambda_t} \sum_{a=0}^{\nu_t-1} \frac{\hat{\mathbf{J}}_t^a}{a!}\log^a\left(\frac{y}{x}\right),$$
where
$$\hat{\mathbf{J}}_t = \begin{pmatrix} 0 & 1 & 0 & \ldots & 0 \\ 0 & 0 & 1 & \ldots & 0 \\ \vdots & \vdots & \vdots & \ddots & \vdots \\ 0 & 0 & 0 & \ldots & 1 \\ 0 & 0 & 0 & \ldots & 0 \end{pmatrix}, \quad \hat{\mathbf{J}}_t^2 = \begin{pmatrix} 0 & 0 & 1 & \ldots & 0 \\ 0 & 0 & 0 & \ldots & 0 \\ \vdots & \vdots & \vdots & \ddots & \vdots \\ 0 & 0 & 0 & \ldots & 0 \\ 0 & 0 & 0 & \ldots & 0 \end{pmatrix}, \ldots.$$

So, $\mathbf{\Sigma}_{11} = \mathbf{0}$. For $g, h = 2, \ldots, s$, $\mathbf{\Sigma}_{1g}^* = \mathbf{\Sigma}_{g1} = \mathbf{0}$,
$$\mathbf{\Sigma}_{gh} = \int_0^1 dx \left[\int_x^1 \frac{1}{y}\left(\frac{y}{x}\right)^{\mathbf{J}_g^*} dy\right]\mathbf{T}_g^*\mathbf{\Sigma}_2\mathbf{T}_h\left[\int_x^1 \frac{1}{y}\left(\frac{y}{x}\right)^{\mathbf{J}_g} dy\right],$$
and the $(a,b)$-element of $\mathbf{\Sigma}_{gh}$ is
$$\sum_{a'=0}^{a-1}\sum_{b'=0}^{b-1}[\mathbf{T}_g^*\mathbf{\Sigma}_i\mathbf{T}_h]_{a-a',b-b'}$$
$$\times \int_0^1 \left[\int_x^1 \frac{(y/x)^{\bar{\lambda}_g}}{a'!y}\log^{a'}\left(\frac{y}{x}\right)dy\right]\left[\int_x^1 \frac{(y/x)^{\lambda_h}}{b'!y}\log^{b'}\left(\frac{y}{x}\right)dy\right]dx$$
$$= \sum_{a'=0}^{a-1}\sum_{b'=0}^{b-1}[\mathbf{T}_g^*\mathbf{\Sigma}_2\mathbf{T}_h]_{a-a',b-b'}$$
$$\times \left[\sum_{l=0}^{a'}\frac{(b'+l)!}{l!b'!}(1-\bar{\lambda}_g)^{-(a'-l+1)}(1-\bar{\lambda}_g-\lambda_h)^{-(b'+l+1)}\right.$$
$$\left.+ \sum_{l=0}^{b'}\frac{(a'+l)!}{l!a'!}(1-\bar{\lambda}_g)^{-(b'-l+1)}(1-\bar{\lambda}_g-\lambda_h)^{-(a'+l+1)}\right],$$
where $[\mathbf{T}_g^*\mathbf{\Sigma}_2\mathbf{T}_h]_{a',b'}$ is the $(a',b')$-element of the matrix $[\mathbf{T}_g^*\mathbf{\Sigma}_2\mathbf{T}_h]$.

# REFERENCES

[1] ANDERSEN, J., FARIES, D. and TAMURA, R. N. (1994). Randomized play-the-winner design for multi-arm clinical trials. *Comm. Statist. Theory Methods* **23** 309–323.

[2] ATHREYA, K. B. and KARLIN, S. (1967). Limit theorems for the split times of branching processes. *J. Math. Mech.* **17** 257–277. MR216592

[3] ATHREYA, K. B. and KARLIN, S. (1968). Embedding of urn schemes into continuous time Markov branching processes and related limit theorems. *Ann. Math. Statist.* **39** 1801–1817. MR232455

[4] BAI, Z. D. and HU, F. (1999). Asymptotic theorem for urn models with nonhomogeneous generating matrices. *Stoch. Process. Appl.* **80** 87–101. MR1670107

[5] BAI, Z. D. and HU, F. (2005). Asymptotics in randomized urn models. *Ann. Appl. Probab.* **15** 914–940. MR2114994

[6] BAI, Z. D., HU, F. and SHEN, L. (2002). An adaptive design for multi-arm clinical trials. *J. Multivariate Anal.* **81** 1–18. MR1901202

[7] BAI, Z. D., HU, F. and ZHANG, L. X. (2002). The Gaussian approximation theorems for urn models and their applications. *Ann. Appl. Probab.* **12** 1149–1173. MR1936587

[8] BEGGS, A. W. (2005). On the convergence of reinforcement learning. *J. Econom. Theory* **122** 1–36. MR2131871

[9] BENAÏM, M., SCHREIBER, S. J. and TARRÉS, P. (2004). Generalized urn models of evolutionary processes. *Ann. Appl. Probab.* **14** 1455–1478. MR2071430

[10] EGGENBERGER, F. and PÓLYA, G. (1923). Über die statistik verketetter vorgänge. *Z. Ang. Math. Mech.* **1** 279–289.

[11] EISELE, J. (1994). The doubly adaptive biased coin design for sequential clinical trials. *J. Statist. Plann. Inference* **38** 249–262. MR1256599

[12] EISELE, J. and WOODROOFE, M. (1995). Central limit theorems for doubly adaptive biased coin designs. *Ann. Statist.* **23** 234–254. MR1331666

[13] EREV, I. and ROTH, A. (1998). Predicting how people play games: Reinforcement learning in experimental games with unique, mixed strategy equilibria. *Amer. Econ. Rev.* **88** 848–881.

[14] HALL, P. and HEYDE, C. C. (1980). *Martingale Limit Theory and Its Applications*. Academic Press, London. MR624435

[15] HU, F. and ROSENBERGER, W. F. (2003). Optimality, variability, power: Evaluating response-adaptive randomization procedures for treatment comparisons. *J. Amer. Statist. Assoc.* **98** 671–678. MR2011680

[16] HU, F. and ZHANG, L.-X. (2004). Asymptotic properties of doubly adaptive biased coin designs for multi-treatment clinical trials. *Ann. Statist.* **32** 268–301. MR2051008

[17] JANSON, S. (2004). Functional limit theorems for multitype branching processes and generalized Pólya urns. *Stoch. Process. Appl.* **110** 177–245. MR2040966

[18] JENNISON, C. and TURNBULL, B. W. (2000). *Group Sequential Methods with Applications to Clinical Trials*. Chapman and Hall/CRC Press, Florida. MR1710781

[19] JOHNSON, N. L. and KOTZ, S. (1977). *Urn Models and Their Applications*. Wiley, New York. MR488211

[20] KOTZ, S. and BALAKRISHNAN, N. (1997). Advances in urn models during the past two decades. In *Advances in Combinatorial Methods and Applications to Probability and Statistics* (N. Balakrishnan, ed.). Birkhäuser, Boston. MR1456736

[21] MARTIN, C. F. and HO, Y. C. (2002). Value of information in the Pólya urn process. *Inform. Sci.* **147** 65–90. MR1940746

[22] MELFI, V. F. and PAGE, C. (1998). Variability in adaptive designs for estimation of success probabilities. In *New Developments and Applications in Experimental Design* (N. Flournoy, W. F. Rosenberger and W. K. Wong, eds.) 106–114. IMS, Hayward, CA. MR1835121
[23] MELFI, V. F. and PAGE, C. (2000). Estimation after adaptive allocation. *J. Statist. Plann. Inference* **87** 353–363. MR1771124
[24] ROBBINS, H. (1952). Some aspects of the sequential design of experiments. *Bull. Amer. Math. Soc.* **58** 527–535. MR50246
[25] ROSENBERGER, W. F. and LACHIN, J. M. (2002). *Randomization in Clinical Trials: Theory and Practice.* Wiley, New York. MR1914364
[26] ROSENBERGER, W. F., STALLARD, N., IVANOVA, A., HARPER, C. N. and RICKS, M. L. (2001). Optimal adaptive designs for binary response trials. *Biometrics* **57** 909–913. MR1863454
[27] SMYTHE, R. T. (1996). Central limit theorems for urn models. *Stoch. Process. Appl.* **65** 115–137. MR1422883
[28] SMYTHE, R. T. and ROSENBERGER, W. F. (1995). Play-the-winner designs, generalized Pólya urns, and Markov branching processes. In *Adaptive Designs* (N. Flournoy and W. F. Rosenberger, eds.) 13–22. IMS, Hayward, CA. MR1477669
[29] THOMPSON, W. R. (1933). On the likelihood that one unknown probability exceeds another in view of the evidence of the two samples. *Biometrika* **25** 285–294.
[30] WEI, L. J. (1979). The generalized Pólya's urn design for sequential medical trials. *Ann. Statist.* **7** 291–296.
[31] WEI, L. J. and DURHAM, S. (1978). The randomized pay-the-winner rule in medical trials. *J. Amer. Statist. Assoc.* **73** 840–843.
[32] WINDRUM, P. (2004). Leveraging technological externalities in complex technologies: Microsoft's exploitation of standards, in the browser wars. *Research Policy* **33** 385–394.

L.-X. ZHANG
DEPARTMENT OF MATHEMATICS
ZHEJIANG UNIVERSITY
HANGZHOU 310028
PEOPLE'S REPUBLIC OF CHINA
E-MAIL: lxzhang@mail.hz.zj.cn

F. HU
DEPARTMENT OF STATISTICS
UNIVERSITY OF VIRGINIA
HALSEY HALL
CHARLOTTESVILLE, VIRGINIA 22904-4135
USA
E-MAIL: fh6e@virginia.edu

S. H. CHEUNG
DEPARTMENT OF STATISTICS
CHINESE UNIVERSITY OF HONG KONG
SHATIN, N.T.
HONG KONG
PEOPLE'S REPUBLIC OF CHINA
E-MAIL: shcheung@sta.cuhk.edu.hk